\documentstyle{article}
\begin{document}

\newtheorem{lemma}{Lemma}
\newtheorem{theorem}{Theorem}
\newtheorem{definition}{Definition}
\newtheorem{corollary}{Corollary}
\newtheorem{problem}{Problem}
\newtheorem{conjecture}{Conjecture}
\newtheorem{proposition}{Proposition}
\def\R{{\bf R}}
\def\C{{\bf C}}
\def\Z{{\bf Z}}

\author{Iskander A. Taimanov
\thanks{Institute of Mathematics, 630090 Novosibirsk, Russia;
e-mail: taimanov\@math.nsc.ru}
}
\title
{The Weierstrass representation of spheres in ${\bf R}^3$,\\
the Willmore numbers, and soliton spheres}
\date{}
\maketitle

\begin{center}
{\bf \S 1. Introduction}
\end{center}

In the present article we consider Weierstrass representations of spheres
in $\R^3$. An existence of a global Weierstrass representation 
for any compact oriented surface of genus $g \geq 1$ has been established
in \cite{T3,T4} and this proof, in fact, works for spheres also. Being mostly
interested in relations of these representations to the spectral theory
and in possibilities to apply the spectral theory to differential geometry
(\cite{T1,T2}), we preferred to consider the case of spheres separately 
because in this case the spectral theory of Dirac operators 
\begin{equation}
{\cal D} =
\left(
\begin{array}{cc}
0 & \partial \\
- \bar{\partial} & 0
\end{array}
\right)
+
\left(
\begin{array}{cc}
U & 0 \\
0 & U
\end{array}
\right)
\label{1}
\end{equation}
is more developed.

In \S 2 we consider two different representations of spheres: as a plane
completed by adding a point at infinity (a plane representation) and
as a cylinder completed by adding a couple of ``infinities'' 
(a cylindric representation). For both of them we describe the data of 
Weierstrass representations which are the potential of a representation
and a ``wave function'' $\psi$ satisfying the equation
\begin{equation}
{\cal D} \psi = 0
\label{2}
\end{equation}
and some analytic conditions.
These data of spectral theory origin are in one-to-one
correspondence with immersed spheres in $\R^3$ (Theorems 1-4). 
This gives a straightforward procedure for constructing immersions in terms of 
zero-eigenfunctions of Dirac operators on a plane and  on an infinite 
two-dimensional cylinder. 

In \S 3 we consider spheres with one-dimensional potentials which
means that in some cylindric representation the potential $U$ of ${\cal D}$
depends on one variable. We prove that a sphere of revolution is uniquely
reconstructed from the potential only (Theorem 5), 
describe all spheres with one-dimensional potentials in terms of 
the Jost functions (Theorem 6), and prove that
\begin{equation}
{\cal W}(\Sigma) \geq 4\pi \left(\dim_{\bf H} \mbox{Ker}\, {\cal D} \right)^2
\label{3}
\end{equation}
where ${\cal W}$ is the Willmore functional and ${\cal D}$ is a Dirac
operator acting on a spinor bundle over a sphere $\Sigma$ (Theorem 7). 
We conjecture that this estimate is valid for all spheres.

In \S 4 we consider a special class of spheres with one-dimensional 
potentials -- spheres with soliton (or reflectionless)
potentials. This integrable case gives many interesting examples, and,
in particular, for each $N$ an equality in (\ref{3}) is achieved
exactly at some special soliton spheres, the Dirac spheres (\cite{R}).

In Appendix A we give a criterion
distinguishing immersions, of universal coverings of compact surfaces of
higher genera, 
converted into immersions of compact surfaces.

Since in \S 3 and \S 4 we extensively use methods of the inverse 
scattering problem for one-dimensional Dirac operators on a line and 
this subject is not well familiar for geometers, Appendix B contains
a brief exposition of this problem and especially some facts which we
need.

We dedicate the present article to our teacher, S. P. Novikov, on his
60th birthday.

\begin{center}
{\bf \S 2. The Weierstrass representation of spheres}
\end{center}

{\bf 2.1. The local Weierstrass representation.}

First, recall the local Weierstrass representation.
It is based on the following two facts:

{\bf Lemma A.}
(Eisenhart (\cite{Eis}, see also comments in \cite{T3}))
{\sl Let $W$ be a simply connected domain in $\C$, $z_0 \in W$, and
let a vector function $\psi = (\psi_1, \psi_2): W \rightarrow \C^2$ satisfy
(\ref{2}) where ${\cal D}$ is of the form (\ref{1}) and its 
potential $U(z,\bar{z})$ is real-valued. Then the following formulas
$$
X^1(z,\bar{z}) = \frac{i}{2} \int_{z_0}^z
\left( (\bar{\psi}_2^2 + \psi_1^2) d z' -
(\bar{\psi}_1^2 + \psi_2^2) d\bar{z}' \right),
$$
\begin{equation}
X^2(z,\bar{z}) = \frac{1}{2}\int_{z_0}^z
\left(
(\bar{\psi}_2^2 - \psi_1^2) d z' -
(\bar{\psi}_1^2 - \psi_2^2) d\bar{z}'
\right),
\label{4}
\end{equation}
$$
X^3(z,\bar{z}) = \int_{z_0}^z
(\psi_1 \bar{\psi}_2 d z' + \bar{\psi}_1 \psi_2 d\bar{z}')
$$
define an immersion of $W$  into $\R^3$ with
the induced metric 
$D(z,\bar{z})^2 dz d\bar{z}$ of the form
$$
D(z,\bar{z}) = |\psi_1(z,\bar{z})|^2 + |\psi_2(z,\bar{z})|^2
$$
and the Gauss curvature and the mean curvature are}
$$
K(z,\bar{z}) = -\frac{4}{D(z,\bar{z})^2}
\partial \bar{\partial}
\log D(z,\bar{z}) 
\ \ \mbox{and} \ \
H(z,\bar{z}) = 2 \frac{U(z,\bar{z})}{D(z,\bar{z})}.
$$

{\bf Lemma B.} (\cite{T1}) 
{\sl Let $W$ be a domain in $\C$ and let $X: W \rightarrow \R^3$ be a
conformal immersion of $W$ into $\R^3$:
$z \rightarrow X(z,\bar{z}) = (X^1(z,\bar{z}), X^2(z,\bar{z}),
X^3(z,\bar{z}))$. Assume that 
$\partial X^3/\partial z \neq 0$ near $z_0 \in W$. 
Then near $z_0$ the functions
\begin{equation}
\psi_1(z,\bar{z}) =
\sqrt{-\partial \Phi(z,\bar{z})},\ \ \
\psi_2(z,\bar{z}) =
\sqrt{\bar{\partial} \Phi(z,\bar{z})},
\label{5}
\end{equation}
with
$$
\Phi(z,\bar{z}) = X^2(z,\bar{z}) + i X^1(z,\bar{z}),
$$
satisfy (\ref{2})
with $U(z,\bar{z}) = H(z,\bar{z}) D(z,\bar{z})/2$,
where $H$ is the mean curvature and $D^2 dz d\bar{z}$ is the metric 
of the surface $X(W) \subset \R^3$.}

A globalization of this representation requires introducing
spinor bundles, generated by $\psi$, over closed oriented 
surfaces and considering the operator (\ref{1}) as acting on them.
This had been shown in \cite{T1} and an existence of a global Weierstrass
representation has been proved for any $C^3$-regular
compact oriented surface of genus $g \geq 1$ (see Theorem 2 in 
\cite{T4}). The proof of it consists of continuing the sections
(\ref{5}) over the whole spinor bundle and uses the following lemma
whose proof is contained in the proof of Theorem 2 from \cite{T4}:

\begin{lemma}
Let $z_0 \in W$ be a nondegenerate critical point of
the function $X^3$ defined on a domain $W \subset \C$ conformally
immersed into $\R^3$. Then near $z_0$ the branches of (\ref{5}) 
are correctly defined as one-valued functions and do not ramify at $z_0$.
\end{lemma}

It is clear that if 
$$
\psi = 
\left(
\begin{array}{c}
\psi_1 \\ \psi_2
\end{array}
\right)
:W \rightarrow \C
$$ 
satisfies (\ref{2}), then
\begin{equation}
\psi^* = 
\left(
\begin{array}{c}
\bar{\psi}_2 \\ -\bar{\psi}_1
\end{array}
\right)
\label{6}
\end{equation}
also satisfies (\ref{2}). Hence for any
$\lambda, \mu \in \C$, such that $|\lambda|^2 + |\mu|^2 \neq 0$,
the vector function
$$
\Psi_{\lambda,\mu} =
\lambda \psi + \mu \psi^* =
\lambda
\left(
\begin{array}{c}
\psi_1 \\ \psi_2 
\end{array}
\right) + 
\mu
\left(
\begin{array}{c}
\bar{\psi}_2 \\ -\bar{\psi}_1 
\end{array}
\right): W \rightarrow \C 
$$
satisfies (\ref{2}).
Consider all immersions 
$$
X_{\lambda,\mu}: W \rightarrow \R^3
$$ 
given via (\ref{4}) by $\Psi_{\lambda,\mu}$ 
and normalize them by the condition
$$
X_{\lambda,\mu}(z_0) = 0
$$
where a point $z_0$ is fixed.
Denote $\partial X^j_{1,0}/\partial z$ by $N^j$.
It is shown by straightforward computations that
$$
\frac{\partial X^1_{\lambda,\mu}}{\partial z} = 
\frac{(\lambda^2 +\bar{\lambda}^2 + \mu^2 + \bar{\mu}^2)}{2} \cdot N^1 +
\frac{i(-\lambda^2 +\bar{\lambda}^2 + \mu^2 -\bar{\mu}^2)}{2} \cdot N^2 +
i(\bar{\lambda}\bar{\mu} - \lambda\mu) \cdot N^3,
$$
\begin{equation}
\frac{\partial X^2_{\lambda,\mu}}{\partial z} = 
\frac{i(\lambda^2 - \bar{\lambda}^2 + \mu^2 - \bar{\mu}^2)}{2} \cdot N^1 +
\frac{(\lambda^2 + \bar{\lambda}^2 - \mu^2 -\bar{\mu}^2)}{2} \cdot N^2 +
(\bar{\lambda}\bar{\mu} + \lambda\mu) \cdot N^3,
\label{7}
\end{equation}
$$
\frac{\partial X^3_{\lambda,\mu}}{\partial z} = 
i(\mu\bar{\lambda} - \lambda\bar{\mu}) \cdot N^1 +
(-\lambda\bar{\mu}-\mu\bar{\lambda}) \cdot N^2 +
(\lambda\bar{\lambda} - \mu\bar{\mu}) \cdot N^3.
$$

From (\ref{7}) we derive

\begin{lemma}  
1) The transformation $\Psi_{1,0} \rightarrow \Psi_{r,0}$, with $r \in \R$,
generates the homo\-thety
$$
(X^1,X^2,X^3) \rightarrow (r^2X^1,r^2X^2,r^2X^3)
$$
of the immersed surface ;

2) For $|\lambda|^2+|\mu|^2 = 1$ the formulas (\ref{7}) define 
the isomorphism
$$
\rho: \{ |\lambda|^2+|\mu|^2 = 1, \lambda,\mu \in \C\}/\{\pm 1\} \rightarrow
SO(3)
$$
and the immersion $X_{\lambda,\mu}$ is a transformation of
$X_{1,0}$ by the rotation $\rho(\lambda,\mu)$.

For instance, 
$$
\rho(\cos \varphi, \sin \varphi) =
\left(
\begin{array}{ccc}
1 & 0 & 0 \\
0 & \cos 2\varphi & \sin 2\varphi \\
0 & -\sin 2\varphi & \cos \varphi
\end{array}
\right),
$$
$$
\rho(e^{i\theta}, 0) =
\left(
\begin{array}{ccc}
\cos 2\theta &  \sin 2\theta & 0 \\
-\sin 2\theta & \cos 2 \theta & 0 \\
0 & 0 & 1
\end{array}
\right),
$$
$$
\rho(\cos \tau, i\sin \tau) =
\left(
\begin{array}{ccc}
\cos 2\tau &  0 & \sin 2\tau  \\
0 & 1 & 0 \\
-\sin 2\tau & 0 & \cos 2\tau 
\end{array}
\right).
$$
\end{lemma}

{\bf 2.2. A plane representation.}

For constructing a global Weierstrass
representation of a two-sphere immersed into $\R^3$ we may
consider a sphere as a plane completed by a point at infinity or
as an infinite cylinder completed by two infinities.
We analyze both possibilities and start with a {\sl plane representation}.

Let $\Sigma$ be a $2$-sphere immersed into $\R^3$.

Fix a pair of points $\infty_{\pm}$ on $\Sigma$ and define a pair
of charts with conformal parameters $z$, on 
$\C \approx \Sigma \setminus \infty_+$, and
$$
u = -\frac{1}{z},
$$
on $\C \approx \Sigma \setminus \infty_-$,
such that $z(\infty_-) = 0$ and $u(\infty_+) = 0$. 
We have
$$
d u = \frac{d z}{z^2}, \ \ \
\frac{\partial}{\partial u} = z^2 \frac{\partial}{\partial z}.
$$
Introduce also the following functions
$$
\tilde{\psi}_1(u,\bar{u}) =
\sqrt{-\partial_u (X^2 + i X^1)},\ \ \
\tilde{\psi}_2(u,\bar{u}) =
\sqrt{\partial_{\bar{u}} (X^2 +i X^1)}.
$$
Now we arrive at the definition

\begin{definition}
A sphere $\Sigma$, immersed into $\R^3$, possesses a global
(plane) Weierstrass representation if there exist  real
potentials $U(z,\bar{z})$ and $\tilde{U}(u,\bar{u})$
and the vector functions $\psi(z,\bar{z})$ and $\tilde{\psi}(u,\bar{u})$
defined on a covering of $\Sigma$ by a pair of charts with parameters
$z$ and $u=z^{-1}$  such that

1)
\begin{equation}
\left\{
\begin{array}{l}
\tilde{U}(u,\bar{u}) = |z|^2 \cdot U(z,\bar{z}),\\
\tilde{\psi}_1 (u,\bar{u}) = z \cdot \psi_1(z,\bar{z}),\\ 
\tilde{\psi}_2 (u,\bar{u}) = \bar{z} \cdot \psi_2(z,\bar{z})
\end{array}
\right.
\label{8}
\end{equation}
for $z = -1/u \in \C$;

2) the vector functions $\psi$ and $\tilde{\psi}$ satisfy (\ref{2}) 
for the Dirac operators with corresponding potentials $U$ and $\tilde{U}$ and  
for a suitable choice of coordinates in $\R^3$ define by (\ref{4})
an immersion of $\Sigma$.
\end{definition}

Consider the analytic conditions
met by the Dirac operator and the spinor sections $\psi$ (\ref{8}) 
corresponding to an immersion of a two-sphere into $\R^3$.

First, notice that
$\tilde{D}(u,\bar{u})^2 du\, d\bar{u} = D(z,\bar{z})^2|z|^4 du\, d\bar{u}$
near $\infty_+$ which implies that
$$
D(z,\bar{z}) = \frac{C}{|z|^2} + O\left(\frac{1}{|z|^3}\right),
\ \mbox{with $C = \mbox{const} \neq 0$}.
$$
and 
\begin{equation}
U(z,\bar{z}) = \frac{U_+}{|z|^2} + 
O\left(\frac{1}{|z|^3}\right)
\ \mbox{with $U_+ = \mbox{const}$}
\label{9}
\end{equation}
as $z \rightarrow \infty$.
Therefore we conclude
\begin{equation}
|\psi(z,\bar{z})|^2 = 
|\psi_1(z,\bar{z})|^2 + |\psi_2(z,\bar{z})|^2
 = O\left(\frac{1}{|z|^2}\right).
\label{10}
\end{equation}
In fact, taking into account that the point $\infty_+ \in \Sigma$
is regular the last equality is refined as follows
\begin{equation}
|\psi(z,\bar{z}|^2 =   
\frac{C_+}{|z|^2} + O\left(\frac{1}{|z|^3}\right), \ \
\mbox{with $C_+ \neq 0$, as $z \rightarrow \infty$.}
\label{11}
\end{equation} 

Assume that $\Sigma$ is $C^3$-regularly immersed into $\R^3$.
By using the general position argument, we choose coordinates in $\R^3$ such
that all critical points of the function $X^3$ defined on $\Sigma$ are
nondegenerate. Take a point $\infty_+ \in \Sigma$ and introduce
a conformal parameter $z$ on $\C \approx \Sigma \setminus \infty_+$. 
By Lemma 1, the branches of (\ref{5}) do not ramify anywhere and
are correctly defined on $\C$. Now from Lemma B and preceding conversations
imply

\begin{theorem}
Every $C^3$-regular two-sphere $\Sigma$ immersed into $\R^3$
possesses a global (``plane'') Weierstrass representation and
the functions $U(z,\bar{z})$ and $\psi(z,\bar{z})$ defined on 
$\Sigma \setminus \infty_+$ satisfy (\ref{9}) and (\ref{11}).
\end{theorem}

In fact, the conditions (\ref{2}), (\ref{9}), and (\ref{11}) 
distinguish the data of Weierstrass representations of spheres.

\begin{theorem}
Let $U(z,\bar{z})$ be a continuous function and
satisfy (\ref{9}) and let $\psi$ be a solution to (\ref{2})
such that $(|\psi_1|^2 + |\psi_2|^2)$ vanishes nowhere on $\C$ 
and (\ref{11}) holds.
Then $\psi$ defines via (\ref{4}) an immersion of
$\C$ completed to a $C^2$-regular immersion of
$S^2$ into $\R^3$. Moreover $U$ and $\psi$ and the functions
$\tilde{U}$ and $\tilde{\psi}$ constructed from them by (\ref{8}) form
the data of a plane Weierstrass representation of the immersed sphere.
\end{theorem}
 
Proof of Theorem 2.

By the definition of the local Weierstrass representation,
$\psi$ defines up to translations an immersion of
$X: \C \rightarrow \R^3$. Normalize this immersion by
$X(i) = 0 \in \R^3$.

Construct the functions $\tilde{U}$ and $\tilde{\psi}$ from $U$ and $\psi$
by (\ref{8}) and notice that they satisfy the following equation
$$
\left[
\left(
\begin{array}{cc}
0 & \partial_u \\
- \partial_{\bar{u}} & 0
\end{array}
\right)
+
\left(
\begin{array}{cc}
\tilde{U} & 0 \\
0 & \tilde{U}
\end{array}
\right)
\right]
\left(
\begin{array}{c}
\tilde{\psi}_1 \\ \tilde{\psi}_2
\end{array}
\right) = 0.
$$
It follows from (\ref{8}) and (\ref{9}) that $\tilde{U}$ is a real-valued
continuous function on the whole complex plane parameterized by $u \in \C$.
Hence $\tilde{\psi}$ also defines up to translations
an immersion of $\tilde{X}:\C \rightarrow \R^3$. Normalize this immersion 
by $\tilde{X}(i) = 0 \in \R^3$.

Since $u=z^{-1}$, we have 
$$
(\bar{\psi}_2^2 + \psi_1^2) d z -
(\bar{\psi}_1^2 + \psi_2^2) d\bar{z} = 
(\bar{\tilde{\psi}}_2^2 + \tilde{\psi}_1^2) d u -
(\bar{\tilde{\psi}}_1^2 + \tilde{\psi}_2^2) d\bar{u}, 
$$
$$
(\bar{\psi}_2^2 - \psi_1^2) d z -
(\bar{\psi}_1^2 - \psi_2^2) d\bar{z}=
(\bar{\tilde{\psi}}_2^2 - \tilde{\psi}_1^2) d u -
(\bar{\tilde{\psi}}_1^2 - \tilde{\psi}_2^2) d\bar{u},
$$
$$
\psi_1 \bar{\psi}_2 d z + \bar{\psi}_1 \psi_2 d\bar{z}=
\tilde{\psi_1} \bar{\tilde{\psi}}_2 d u + 
\bar{\tilde{\psi}}_1 \tilde{\psi}_2 d\bar{u}
$$
on $\C^* = \C \setminus \{0\}$. 
These formulas imply that
$$
X(z) = \tilde{X}\left(-\frac{1}{z}\right) \ \ \mbox{for $z \in \C^*$}.
$$
Hence, $X$ and $\tilde{X}$
coincide  on $\C^*$ and each of them is regularly continued onto
the corresponding ``infinity point'', on $u=0$ and $z = 0$.
By (\ref{11}), $\tilde{X}$ is regular at $u=0$.

This proves the theorem.

A nice feature of this theorem is that the closedness  problem consisting
in distinguishing immersions of planes which are converted into immersions
of compact surfaces reduces for spheres to the conditions
(\ref{9}) and (\ref{11}) which can be easily checked. 
For surfaces of higher genera this problem is more complicated
(see Appendix A).
 
If we have a solution $\psi$ to (\ref{2}) such that $|\psi|$
decays slower than $|z|^{-1}$ as $z \rightarrow \infty$ then nevertheless 
we may construct by using (\ref{4}) an immersion of a sphere into $\R^3$
but it would not be regular and would have a peak singularity at the
``infinity''. If $|\psi|$ decays faster than $|z|^{-1}$ than
we would have a branch point at ``infinity''. This also occurs when
$|\psi(z,\bar{z})|=0$ at $z \in \C$.

Let now admit branch points and consider more general situation.

Denote by ${\cal E}$ a $\C^2$-bundle
$$
\C^2 \rightarrow {\cal E} \rightarrow S^2
$$
whose sections $\psi$ satisfy (\ref{8}). This is a spinor bundle
obtained as a square root of the complexified tangent bundle.   

\begin{theorem}
Let $U$ satisfy (\ref{8}) and (\ref{9}). Then 

1) ${\cal D}$ acts on sections of ${\cal E}$;

2) solutions to (\ref{2}) satisfying
(\ref{10}) are in one-to-one correspondence with zero-eigenfunctions
of ${\cal D}$;

2) the kernel of ${\cal D}$ is finite-dimensional and
moreover it is even-dimensional.
\end{theorem}

The first and second statements are evident. 
Since ${\cal D}$ is elliptic its kernel is 
finite-dimensional. We know that there exists an automorphism $*$, of
the kernel, given by (\ref{6}). Since $(\psi^*)^* = -\psi$, the kernel 
splits into two-dimensional subspaces invariant under $*$ and therefore
$\dim_{\C}\mbox{Ker}\, {\cal D} = 2n$ with $n$ integer
\footnote{F. Pedit and U. Pinkall proposed to treat 
${\cal E}$ as a quaternion vector bundle and to treat $\mbox{Ker}\, {\cal D}$ 
as a vector space over quaternions identifying $*$ with a 
multiplication by ${\bf j} \in {\bf H}$ (\cite{KPP}).
In this case $\dim_{\bf H}\mbox{Ker}\, {\cal D} = n$.}.
Each section $\psi \in\mbox{Ker}\, {\cal D}$ generates via (\ref{4}) 
an immersed sphere which may have branch points.

{\bf 2.3. A cylindric representation.}
 
Starting with a plane representation we also introduce
a {\sl cylindric representation}, i.e., a representation of an immersed
sphere as an immersed cylinder completed by two points.

Put $z' = \log z  = x' +i y'$ and take a cylinder 
${\cal Z} = \C/i\Z$ with a conformal parameter $z'$ given modulo $2\pi i$.

Let $(U,\psi)$ be the data of a plane representation.
We have
$$
d z' = \frac{1}{z} d z, \ \ 
\frac{\partial}{\partial z'} = z\frac{\partial}{\partial z}.
$$
Consider the following functions
\begin{equation}
\hat{D}(z,\bar{z}') = |z| D(z,\bar{z}), \ \
\hat{U}(z,\bar{z}') = |z| U(z,\bar{z})
\label{12}
\end{equation}
and
\begin{equation}
\hat{\psi}_1(z',\bar{z}') = \sqrt{z} \psi_1(z,\bar{z}), \ \
\hat{\psi}_2(z',\bar{z}') = \sqrt{\bar{z}} \psi_2(z,\bar{z}).
\label{13}
\end{equation}
By straightforward computations it is obtained that

\begin{lemma}
1) The functions $\hat{\psi}$ satisfy the equation
\begin{equation}
\left[
\left(
\begin{array}{cc}
0 & \partial_{z'} \\
- \partial_{\bar{z}'} & 0
\end{array}
\right)
+
\left(
\begin{array}{cc}
\hat{U} & 0 \\
0 & \hat{U}
\end{array}
\right)
\right]
\left(
\begin{array}{c}
\hat{\psi}_1 \\ \hat{\psi}_2
\end{array}
\right) = 0.
\label{14}
\end{equation}

2) There are the asymptotics 
\begin{equation}
\hat{U} = \frac{U_{\pm}}{e^{|x'|}} + O\left(\frac{1}{e^{2|x'|}}\right), \ \  
|\hat{\psi}_1|^2 +|\hat{\psi}_2|^2 = \frac{C_{\pm}}{e^{|x'|}}
+ O\left(\frac{1}{e^{2|x'|}}\right),
\label{15}
\end{equation}
with $U_{\pm}$ and $C_{\pm}$ constants such that $C_{\pm} \neq 0$
as $x' \rightarrow \pm \infty$.
\end{lemma}

Now it is clear how to derive from Theorems 1 and 2 the following
result

\begin{theorem}
1) Any vector function $\hat{\psi}$ satisfying 
(\ref{14}) and (\ref{15}) 
defines via (\ref{4}) an immersion of
${\cal Z}$ into $\R^3$, which is completed to a regular immersion
of a two-sphere.

2) For any two-sphere $\Sigma$ which is $C^3$-immersed into $\R^3$ and
any pair of distinct points $\infty_{\pm} \in \Sigma$ there exists
an immersion of a cylinder ${\cal Z} = \Sigma \setminus \{\infty_{\pm} \}$
such that 

a) this immersion is defined, for a suitable choice of coordinates in $\R^3$,
via (\ref{4}) by functions $\hat{\psi}$ and $\hat{U}$ satisfying 
(\ref{14}) and (\ref{15});

b) by adding a point to each end of the cylinder 
this immersion is completed to a regular immersion of $\Sigma$.
\end{theorem}

It follows from Theorem 1 and Theorem 4 that

\begin{corollary}
Every $C^3$-regular two-sphere $\Sigma$ immersed into $\R^3$ possesses a 
cylindric Weierstrass representation.
\end{corollary}

We mention above that $\psi(z,\bar{z})$ are sections of ${\cal E}$.
The formulas 
(\ref{12}) and (\ref{13}) 
show that vector functions $\hat{\psi}$ on ${\cal Z}$ meeting the 
conditions (\ref{15}) and
\begin{equation}
\hat{\psi} (x',y'+ 2\pi) = -\hat{\psi} (x',y')
\label{16}
\end{equation} 
are sections of ${\cal E}$
and these formulas just establish an equivalence between two different
representations  of ${\cal E}$.
Moreover, these formulas also establish
the equivalence between Dirac operators and, therefore, we have

\begin{corollary}
Solutions $\hat{\psi}$ to (\ref{14}) satisfying (\ref{16}) and
$$
|\hat{\psi}|^2 = O\left(\frac{1}{e^{|x'|}}\right)
\ \ \ \
\mbox{as $x' \rightarrow \pm \infty$}
$$
form the kernel of ${\cal D}$. The dimension of $\mbox{Ker}\, {\cal D}$ 
is finite and even.
\end{corollary} 

We also mention that each section from this kernel generates via (\ref{4}) an 
immersed sphere which may have branch points.

\begin{center}
{\bf \S 3. Spheres with one-dimensional potentials}
\end{center}

{\bf 3.1. The spectral data for one-dimensional potentials.}

In this chapter we consider spheres admitting cylinder representations
with one-dimensional real-valued potentials $U(x)$. This means that being 
defined on an infinite cylinder 
${\cal Z} = \{(x,y): -\infty < x < \infty, 0 \leq y 
\leq 2\pi\}$ the potential $U(z)$ depends on $x$ only.
We assume that $U(x)$ decays exponentially:
$$
U(x) = O\left(\frac{1}{e^{|x|}}\right)
\ \ \ \
\mbox{as $x \rightarrow \pm \infty$.}
$$ 

For $U = 0$  solutions to (\ref{2}) are linear combinations of
$$
\left(
\begin{array}{c} 
0 \\ \exp(ik\bar{z})
\end{array}
\right)
\ \ \ \ \mbox{and} \ \ \ 
\left(
\begin{array}{c} 
\exp(ilz) \\ 0
\end{array}
\right)
$$
with $z = x + iy$ and $k,l \in \C$. Of course, such functions are
defined on ${\cal Z}$ if and only if $k, l \in 2\pi \Z$. 
Hence we look for solutions to (\ref{2}) of the form
\begin{equation}
\psi(x,y) = \exp(k y) \varphi(x,y), \ \ \ \varphi(x,y+2\pi) = \varphi(x,y),
\label{17}
\end{equation}
i.e., we consider solutions which are defined on the universal covering 
$\R^2$ of ${\cal Z}$ and satisfy the periodicity condition
$$
\psi(x, y+2\pi) = \mu \psi(x,y)
$$
with $\mu$ a constant. 
For this ansatz the equation (\ref{2}) reduces to
\begin{equation}
\left[
\left(
\begin{array}{cc}
0 & \partial \\
- \bar{\partial} & 0
\end{array}
\right)
+
\left(
\begin{array}{cc}
U & 0 \\
0 & U
\end{array}
\right)
-
\frac{1}{2}
\left(
\begin{array}{cc}
0 & ik \\
ik & 0
\end{array}
\right)
\right]
\varphi = 0.
\label{18}
\end{equation}
and, decomposing its solutions into Fourier series in $y$ 
$$
\varphi(x,y) = \sum_{m \in \Z} \varphi_m(x)e^{i m y},
$$
we conclude that each $\varphi_m(x)$ satisfies ({\ref{18}) with
$k+im$ substituted for $k$. Hence for studying all solutions to (\ref{2})
of the form (\ref{17}) it is enough to study solutions to
\begin{equation}
\left[
\left(
\begin{array}{cc}
0 & \partial_x \\
- \partial_x & 0
\end{array}
\right)
+
\left(
\begin{array}{cc}
2U & 0 \\
0 & 2U
\end{array}
\right)
-
\left(
\begin{array}{cc}
0 & ik \\
ik & 0
\end{array}
\right)
\right]
\varphi = 0
\label{19}
\end{equation}
depending on $x$ only.

This problem
has been studied in its relation to soliton equations
(see, for instance, \cite{AS,FT,Lamb,TS})
and we give the brief summary of results which we need in the sequel
\footnote{A detailed exposition is given in Appendix B.}.	

Since $U(x)$ decays exponentially as $|x| \rightarrow \infty$,
we have

{\bf Summary.}

{\sl The potential $U(x)$ of an operator
$$
L = 
\left(
\begin{array}{cc}
0 & \partial_x \\
- \partial_x & 0
\end{array}
\right)
+
\left(
\begin{array}{cc}
2U & 0 \\
0 & 2U
\end{array}
\right)
$$
is uniquely reconstructed from the spectral data which are

i) the reflection coefficient $R(k)=b(k)/a(k)$, 
with $k \in \R \setminus \{0\}$,

ii) the poles $\kappa_1,\dots,\kappa_N$ of the transmission coefficient
$T(k)$ with $\mbox{Im}\, k>0$,

iii) some additional quantities
$\lambda_1,\dots,\lambda_N \in \C$, attached to $\kappa_j$.

A reconstruction procedure is given by 
the Marchenko equations (\ref{B16}) and (\ref{B20}).

The poles of $T(k)$ are simple and correspond to
exponentially decaying solutions to (\ref{19}). For each pole $\kappa_j$ 
every such solution is a multiple of
$\varphi_1^+(x,\kappa_j)$ which is a unique solution to the equation
$$
\varphi (x) =
\left(
\begin{array}{c}
0 \\ e^{i\kappa_j x} 
\end{array}\right)
+
\int_x^{+\infty}
\left(
\begin{array}{cc}
0 & -e^{-i \kappa_j (x-x')} \\
e^{i \kappa_j (x-x')} & 0
\end{array}\right)
\cdot 2U(x') \cdot \varphi(x')
\, d x'.
$$
Since $U(x)$ is real-valued, 

a) $\kappa_j$ are symmetric with respect to the imaginagy axis,
if $\kappa_j$ and $\kappa_l = -\overline{\kappa_j}$ are different poles
of $T(k)$ then $\lambda_j = \bar{\lambda}_l$, and if $\mbox{Re}\,\kappa_m = 0$
then $\lambda_m \in \R$ ;

b) $R(k) = \overline{R(-k)}$.}
     
We recall the definition of the Kruskal integrals:
$$
I_n(U) = \int_{-\infty}^{+\infty} U(x)q_n(x)\, d x,
$$
where
$$
q_1(x) = U(x)
$$
and the other quantities $q_j(x)$ are defined by the recursion relation
$$
q_{j+1}(x) = -i\frac{d q_j(x)}{d x} - 4U(x)\sum_{m=1}^{j-1}q_m(x)q_{j-m}(x).
$$
The first of them is obviously the squared $L_2$-norm of $U(x)$:
$$
I_1(U) = \int_{-\infty}^{+\infty} U^2(x)\, d x.
$$

These quantities are related to the spectral data via the trace formulas
(see \cite{FT}, the formulas (7.20) and (7.21) in chapter 1 of Part I
\footnote{In \cite{FT} these formulas are written in terms of $\kappa$, 
$\psi(x)$, and $\lambda$ which, in our notation,  are
$(-4)$, $-iU(x)$, and $-2k$, respectively.}): 
\begin{equation}
I_n(U) = 
-\frac{1}{4\pi} \int_{-\infty}^{+\infty} 
\log(1-|b(k)|^2)(-2k)^{n-1}\,d k + \frac{i2^{n-2}}{n} \sum_{j=1}^N
(\bar{\kappa}_j^n - \kappa_j^n)
\label{20}
\end{equation}
where $b(k)$ is the ratio of $R(k)$ and $T(k)$ and satisfies
the inequality
$$
0 < |b(k)| < 1.
$$
For $n=1$ we have
\begin{equation}
\int_{-\infty}^{+\infty} U^2(x) \, d x = 
-\frac{1}{4\pi} \int_{-\infty}^{+\infty} 
\log(1-|b(k)|^2)\,d k + \sum_{j=1}^N \mbox{Im}\, \kappa_j
\label{21}
\end{equation}
and we conclude that
\begin{equation} 
\int_{-\infty}^{+\infty} U^2(x) \, d x
\geq  \sum_{j=1}^N \mbox{Im}\, \kappa_j
\label{22}
\end{equation}
and an equality in (\ref{22}) is achieved exactly at {\sl reflectionless}
potentials, i.e., $b(k) \equiv 0$ which is equivalent to $R(k) \equiv 0$.

{\bf 3.2. Construction of spheres with one-dimensional potentials.
A reconstruction of a sphere of revolution from its potential.}

By Corollary 1, every sphere $\Sigma$ regularly immersed into
$\R^3$ possesses a cylindric Weierstrass representation. 
In this subchapter we describe spheres which admits cylindric 
representations with potentials depending on $x$ only. 
The simplest and most important examples are spheres of revolution
(\cite{T2}). 

Let $U(x)$ be a potential of an immersed sphere $\Sigma$. By Lemma 4,
it decays exponentially and we may apply the spectral theory of $L$
exposed in 3.1 and Appendix B. In particular, all  
exponentially decaying solutions to (\ref{18}) are linear combinations
of $\varphi^+_1(x,\kappa_j)$ and their $*$-transforms (\ref{6}).
 
It is known that for any sphere of revolution $\Sigma$ 
there exists a cylindric
representation such that its potential is one-dimensional and
$\Sigma$ is immersed into $\R^3$ via (\ref{4}) where
$\psi$ has the form
$$
\psi(x,y) = \varphi(x) e^{i y /2}
$$
(see \cite{T2}). Let $U(x)$ be the potential of this representation  
of $\Sigma$. Then we have
$$
\psi_{\lambda,\mu}(x,y) = \lambda \left(\varphi^+_1(x,i/2)e^{i y/2}\right) + 
\mu \left(\varphi^+_1(x,i/2)e^{i y/2}\right)^*
$$
where $\lambda, \mu \in \C$.
By Lemma 2, for different $\lambda$ and $\mu$ such immersions are transformed  
one into another by homotheties of spheres and their rigid motions in
$\R^3$. Since the potentials of different Weierstrass representations are
reconstructed one from another by the formulas (\ref{8}) and (\ref{12}), we
conclude

\begin{theorem}
Any sphere of revolution without branch points 
is uniquely defined (up to homotheties of the sphere
and rigid motions in $\R^3$) by the potential of any of its Weierstrass
representations.
\end{theorem}

The condition on absence of branch points is added just for 
the following reason.
Notice that we may consider the linear combination
$$
\psi_{\lambda,\mu,\kappa}(x,y) = 
\lambda \left(\varphi^+_1(x,\kappa)e^{\kappa y}\right)
+ \mu \left(\varphi^+_1(x,\kappa)e^{\kappa y}\right)^* 
$$ 
with $\kappa = i n/2$ and $n > 1$.
Then  the sphere constructed from $\psi_{\lambda,\mu,\kappa}$ via 
(\ref{4}) would be an $n$-sheeted covering of a sphere of revolution with 
branch points at the infinities.

For general spheres with one-dimensional potentials 
the statement of Theorem 5 does not valid.

Let $\kappa_1,\dots,\kappa_N$ be the poles,
of the transmission coefficient $T(k)$, coming into the spectral data of
$U(x)$ and divide them into three groups
$$
\kappa_j = \frac{i n_j}{2}\ \  
\mbox{with $n_j$ an odd positive integer for} \ \ 
1 \leq j \leq L,
$$ 
$$
\kappa_j = \frac{i n_j}{2}\ \  
\mbox{with $n_j$ an even positive integer for} \ \ 
L+1 \leq j \leq M,
$$ 
and $\kappa_j$ is not of the form $i n /2$ with $n$ integer for 
$j \geq M+1$.
Put 
$$
\psi_j(x,y) = \varphi^+_1(x,\kappa_j)e^{\kappa_j y}.
$$ 

It is clear that for $j \geq M+1$ the functions $\psi_j$ and
$\psi^*_j$ are no periodic and no antiperiodic in $y$. Therefore
squares of linear combinations of such functions are not defined on ${\cal Z}$ 
and do not generate via (\ref{4}) immersions of cylinders.

By Corollary 2, since  $\psi_j$ and $\psi^*_j$ 
satisfy (\ref{16}) for $j \leq L$, they are 
sections of ${\cal E}$ and we conclude

\begin{lemma}
$\mbox{Ker}\, {\cal D}$ is spanned by 
$\psi_j(x,y)$ and $\psi^*_j(x,y)$ where $j \leq L$.
\end{lemma}

Any linear combination $\psi(x,y)$ of $\psi_j$ and $\psi^*_j$ for
$L+1 \leq j \leq M$ also generate via (\ref{4}) an immersion of a sphere. 
It is easy to see that
if for all $\psi_j$ and $\psi^*_k$ coming into this combination the 
frequencies $n_j$ are represented in the form
$$
n_j = 2^k l_j
$$
with $l_j$ odd integers then the immersion would be a $2^k$-sheeted covering
over its image with branch points of order $2^k$ at infinities.
The potential of the representation of a covered sphere given by the
function $\psi'(x,y) = \psi(x/2^k,y/2^k)$ would be
$U'(x) = U(x/2^k)$. Otherwise the immersion would have branch points of odd 
order at the infinities. 

Since by the definition, a vector function $\psi$ coming into a
cylindric representation belongs to $\mbox{Ker}\, {\cal D}$, we 
summarize these conversations as follows

\begin{theorem}
Let ${\bf a} = (a_1,\dots,a_{2L}) \in \C^{2L}\setminus \{0\}$.
Then the function
\begin{equation}
\psi_{\bf a}(x,y) = a_1 \psi_1(x,y) + \dots + a_L\psi_L(x,y) + 
a_{L+1} \psi^*_1(x,y) + \dots + a_{2M}\psi^*_L(x,y)
\label{23}
\end{equation}
defines via (\ref{4}) an immersed sphere $\Sigma_{\bf a}$ in  
$\R^3$. 

If there exist non-zero coefficients $a_j$ and $a_k$ such that 
$(j-k) \neq \pm L$ then $\Sigma_{\bf a}$ is not a sphere of revolution.

These spheres are exactly the spheres which have $U(x)$ as a potential of
some of their cylindric Weierstrass representations.
\end{theorem}

{\bf 3.3. The Willmore functional via the trace formula and 
the Willmore numbers.}

Consider the following problem

\begin{problem}
How to estimate the dimension of $\mbox{Ker}\, {\cal D}$ ?
\end{problem}

For a Dirac operator with a one-dimensional potentials (in some cylindric
representation of a sphere) the trace formula (\ref{21}) enables us to give
a precise estimate.  

Indeed, by Lemma 4, $\mbox{Ker}\, {\cal D}$ is spanned by
$\varphi^+_1(x)e^{\kappa_j y}$ and 
$(\varphi^+_1(x)e^{\kappa_j y})^*$ with $\kappa_j$ of the form
\begin{equation}
\kappa_j = \frac{i (2n_j+1)}{2}
\label{24}
\end{equation}
with $n_j$ nonnegative integers. 
By (\ref{21}) and (\ref{22}), we have
\begin{equation}
\int_{\infty}^{+\infty} U^2(x) \, d x \geq 
\frac{1}{2}\sum_{j=1}^L (2n_j+1)
\label{25}
\end{equation}
and an equality in (\ref{25}) is achieved exactly in the case 
then $U(x)$ is a reflectionless potential and   
the whole discrete spectrum with $\mbox{Im}\, \kappa > 0$ consists of
eigenvalues of the form (\ref{24}). 
We conclude 

\begin{theorem}
Let ${\cal D}$ be a Dirac operator (\ref{1}) on ${\cal E}$
with a one-dimensional potential $U(x)$ 
in some cylindric representation
(for instance, a Dirac operator generating an immersion of a sphere of
revolution). If
$$
\frac{\dim_{\C}\mbox{Ker}\, {\cal D}}{2} = 
\dim_{\bf H}\mbox{Ker}\, {\cal D} \geq N,
$$
then
\begin{equation}
\int_{-\infty}^{+\infty} U^2(x) \, d x \geq \frac{N^2}{2}.
\label{26}
\end{equation}
\end{theorem}

Indeed, any level $\kappa = i(2n+1)/2$ may be filled just by one
eigenfunction of the form $\varphi^+_1(x,\kappa_j)e^{\kappa_j y}$ 
(see Appendix B) and, given 
$N = \dim_{\bf H}\mbox{Ker}\, {\cal D}$, the left-hand side in (\ref{26})
achieves its minimal possible value if just first $N$ levels are filled.
This means that 
$$
\frac{1}{2}\sum_{j=1}^N (2n_j+1) = 
\frac{1}{2}(1 + 3 + \dots + (2N-1)) = \frac{N^2}{2}.
$$ 
This proves the theorem.

An example of the Dirac spheres, 
constructed by U. Pinkall and J. Richter (\cite{R}),
shows that for each $N$ an equality in (\ref{26}) is achieved at
the potential of such sphere and therefore 
\begin{corollary}
For any $N$ the estimate (\ref{26}) is precise and 
an equality is achieved at 
$$
U_N(x) = \frac{N}{2 \cosh x}.
$$
\end{corollary}
We discuss the spectral data of such potentials in \S 4.

To rewrite (\ref{26}) for general spheres it needs to integrate
the left-hand side over $y$ also and obtain an integral over $S^2$
(by (\ref{8}) and (\ref{12}) this integral is correctly defined for
any representation):
\begin{equation}
\int_{\Sigma} U^2(z,\bar{z}) \,d x \wedge d y \geq \pi N^2. 
\label{27}
\end{equation}

We would like to conjecture that
\begin{conjecture}
The estimate (\ref{27}) holds for any Dirac operator on ${\cal E}$.
\end{conjecture}

In fact, by Theorem 4, the dimension of $\mbox{Ker}\, {\cal D}$ measures
the dimension of a family of isopotential spheres in $\R^3$. 
If there exists $\psi \in \mbox{Ker}\, {\cal D}$ such that $\psi$
vanishes nowhere then it generates an immersion of a sphere without branch
points. It is also known that for any compact surface $\Sigma$
immersed via (\ref{4})
the value of the Willmore functional ${\cal W}$, 
an integral of a squared mean curvature, is given by
\begin{equation}
{\cal W}(\Sigma) = \int H^2\, d\mu = 
4\int_{\Sigma} U^2(z,\bar{z}) \,d x \wedge d y
\label{28}
\end{equation}
(see \cite{T1}), i.e., a multiple of the squared $L_2$-norm of 
the potential $U$. We have
$$
\int_{\Sigma} (H^2-K) \,d\mu = \int_{\Sigma} 
\left(\frac{k_1-k_2}{2}\right)^2 \,d\mu\geq 0
$$
where $k_j$ are the principal curvatures and $K$ is the Gauss curvature of
$\Sigma$. By the Gauss--Bonnet theorem, for spheres
$$
\int_{\Sigma} K \,d\mu = 4\pi
$$
and this implies
$$
\int_{\Sigma} H^2\, d\mu \geq 4\pi. 
$$
Therefore we conclude that
$$
\int_{\Sigma} U^2(z,\bar{z}) \,d x \wedge d y \geq \pi.
$$
But this is just the inequality (\ref{27}) for $N=1$ and 
an equality is achieved at $U_1(x)$ which is the potential of the 
unit sphere (\cite{T2}). Recalling the recent result of F. Pedit
and U. Pinkall obtained by methods of the so-called quaternionic
algebraic geometry (\cite{PP}), we conclude

\begin{proposition}
Assume that there exists a zero-eigenfunction $\psi$ of ${\cal D}$
such that $\psi$ defines via (\ref{4}) a regular immersion of $S^2$ into
$\R^3$. Then Conjecture 1 is valid 

1) for $N=1$ (Gauss--Bonnet);

2) for $N=2$ (Pedit--Pinkall).
\end{proposition}  

For such operators the estimate (\ref{27})
in terms of ${\cal W}$ takes the form
\begin{equation}
{\cal W}(\Sigma) \geq 4\pi N^2.
\label{29}
\end{equation}
We show in \S 4 
the dimension of $\mbox{Ker}\, {\cal D}$ can not be estimated from below
in terms of the Willmore functional.

This treatment of (\ref{27}) fits into a general approach to estimates
for the Willmore functional based on the Weierstrass representation
(see \cite{T4} where the spectral approach 
for the Willmore conjecture for tori is introduced). The right-hand sides of
(\ref{29}), {\sl the Willmore numbers}, measure not only the
existence of a sphere with given value of the Willmore functional  
but the dimension of a family of isopotential spheres.
It looks natural that for given dimension of $\mbox{Ker}\, {\cal D}$ 
the Willmore functional has to attain its minimal possible value on 
a very symmetric operator which has to have a one-dimensional potential,
i.e., to be of the form covered by Theorem 7. The analogous idea is
discussed in \cite{T4} for the Willmore conjecture for tori.

\begin{center}
{\bf \S 4. Soliton spheres}
\end{center}

{\bf 4.1. Solving the Marchenko equations for reflectionless potentials.}

We call a sphere reflectionless if it admits a Weierstrass representation
with a one-dimensional reflectionless potential $U(x)$. This means that
the reflection coefficient $R(k)$ of $U(x)$ vanishes identically, i.e.,
$R(k) \equiv 0$, and the spectral data of $U(x)$ are just

1) a half, of a discrete spectrum of $L$
with $U(x)$ its potential, lying in the upper-half plane:
$\kappa_1,\dots,\kappa_N$; this spectrum is symmetric with respect to
the imaginary axis;

2) some quantities $\lambda_j$ corresponding to $\kappa_j$ such that if 
$\mbox{Re}\, \kappa_j = 0$ then $\lambda_j \in \R$ and if
$\kappa_k$ and $\kappa_l=-\bar{\kappa}_k$ do not coincide 
then $\lambda_k$ = $\bar{\lambda}_l$.

The potential $U(x)$ is reconstructed from the spectral data via
the Marchenko equations (see Appendix B). 
For reflectionless potentials solutions to these equations can be 
found explicitly. We explain this procedure following \cite{Lamb}. 

Given the spectral data for a reflectionless potential, 
consider the following ansatz:
$$
B_j(x,y) = \langle B_j(x) | T(y) \rangle
$$ 
where
$$
T(z) = (e^{i\kappa_1 z},\dots,e^{i\kappa_N z})
$$
and
$$
\langle u | v \rangle = u_1 v_1 + \dots + u_N v_N
$$
is a standard inner product.
Then $\Omega$ takes the form
$$
\Omega(x+y) = \langle \Psi(x) | T(y) \rangle  
$$
where
$$
\Psi(z) = (-\lambda_1 e^{i\kappa_1 z}, \dots, -\lambda_N e^{i\kappa_N z}).
$$
In terms of these functions the Marchenko equations (\ref{B16})
are written as
$$
\langle B_2(x) | T(y) \rangle + 
\int_x^{+\infty} \langle B_1(x) | T(x') \rangle
\langle \Psi(x') | T(y) \rangle\, d x' = 0
$$
and
$$
\langle \Psi(x) | T(y) \rangle - \langle B_1(x) | T(y) \rangle  +  
\int_x^{+\infty} \langle B_2(x) | T(x') \rangle
\langle \Psi(x') | T(y) \rangle\, d x' = 0.
$$
Introduce the matrix
$$
M(x) = \int_x^{+\infty} | T(x') \rangle \langle \Psi(x')|\, d x'
$$
(here we use Dirac's notation  treating the inner product 
as a product of a bra vector $\langle u |$ and a ket vector $| v \rangle$,
\cite{Dirac}), rewrite the Marchenko equations as
$$
\langle B_2(x) + B_1(x)M(x) | T(y) \rangle = 
\langle  \Psi(x) - B_1(x) + B_2(x)M(x) | T(y) \rangle = 0
$$
and finally arrive at the following form of them
\begin{equation}
B_2(x) + B_1(x)M(x) =  \Psi(x) - B_1(x) + B_2(x) M(x) = 0. 
\label{30}
\end{equation}
The entries of $M(x)$ are simply computed
$$
M_{j k}(x) = 
\int_x^{+\infty} T_j(x')\Psi_k(x') \, d x' = 
\frac{\lambda_k}{i(\kappa_j +\kappa_k)}e^{i(\kappa_j+\kappa_k)x}
$$
and (\ref{30}) implies that
$$
B_1(x) = \Psi(x) \cdot (1+M^2(x))^{-1}
$$
and
$$
B_2(x) = - B_1(x) \cdot M(x).
$$
By (\ref{B19}), we derive
$$
U(x) = - \langle  \Psi(x) \cdot (1+M^2(x))^{-1} | T(x) \rangle.
$$
Now represent $B_1(x,x)$ and $B_2(x,x)$ as follows
$$
B_1(x,x) = 
- \mbox{Tr}\,\left[ \frac{d M(x)}{d x} \cdot (1 + M^2(x))^{-1} \right],
$$
$$
B_2(x,x) = 
\mbox{Tr}\,\left[ \frac{d M(x)}{d x} \cdot M(x) \cdot 
(1 + M^2(x))^{-1} \right].
$$
It follows from (\ref{B19}) and (\ref{B21}) that
$$
2U^2(x) + i\frac{d U(x)}{d x} = 
\frac{d}{d x} \mbox{Tr}\,
\left[ \frac{d M(x)}{d x}(M(x) +i) (1 + M^2(x))^{-1}\right],
$$
and, since $1 + M^2(x) = (1+ i M(x))(1 - i M(x))$, we have
$$
2U^2(x) + i\frac{d U(x)}{d x} = 
\frac{d}{d x}\mbox{Tr}\, 
\left[ \frac{d (1 +i M(x))}{d x} (1 + i M(x))^{-1} \right].
$$
Using the well-known identity
$$
\frac{d}{d x}\, \log \det A(x) = \mbox{Tr}\, 
\left( \frac{d A(x)}{d x} \cdot A^{-1}(x) \right),
$$
we obtain
$$
2U^2(x) + i\frac{d U(x)}{d x} = 
\frac{d^2}{d x^2} \log \det (1 + i M(x)).
$$
Since $U(x)$ is real-valued and fast decaying, we have
$$
U(x) = \frac{d}{d x} \mbox{Im}\, \log \det (1 + i M(x)) =
\frac{d}{d x} \arctan 
\frac{\mbox{Im}\, \det (1 + i M(x))}{\mbox{Re}\, \det (1 + i M(x))}.
$$
For reflectionless potentials $\varphi^+_1(x,k)$
is simply written. Put
$$
W(x,k) = \int_x^{+\infty} T(x') e^{i k x'} \, d x' = 
\left(
\frac{i}{\kappa_1 + k}e^{i ( \kappa_1 + k)x}, \dots , 
\frac{i}{\kappa_N + k}e^{i ( \kappa_N + k)x}
\right).
$$
and, by (\ref{B11}), obtain
\begin{equation}
\varphi^+_1(x,k) = 
\left(
\begin{array}{c}
\langle \Psi(x) \cdot (1+M^2(x))^{-1} | W(x,k) \rangle  \\ 
e^{i k x}  -
\langle \Psi(x) \cdot (1+M^2(x))^{-1} M(x) | W(x,k) \rangle
\end{array}
\right).
\label{31}
\end{equation}

{\bf 4.2. Construction and properties of reflectionless spheres.}

We consider some explicit examples of reflectionless spheres
which are spheres with $N$-soliton potentials.

{\bf 4.2.1. $N=1,\, \kappa_1 = \frac{i}{2}$.}

In this case 
$$
M(x) = -\lambda e^{-x}
$$
where $\lambda = \lambda_1 \in \R \setminus \{ 0 \}$
and we obtain
$$
U(x,\lambda) = \frac{\lambda e^{-x}}{1 + \lambda^2 e^{-2x}}.
$$
Since for an immersed surface $U(x)$ is defined up
to a sign, we assume that $\lambda = e^{-a} > 0$ and
derive $U(x,\lambda) = U_1(x+a)$ with
$$
U_1(x) = \frac{1}{2 \cosh x}.
$$
This is the potential of a round sphere and it is easily checked
that $\varphi^+_1(x,i/2)e^{i y/2}$ defines it via (\ref{4}).

{\bf 4.2.2. $N = 2 ,\, \kappa_1 = \frac{i}{2},\, 
\kappa_2 = \frac{3i}{2}$.}

A general potential corresponding to this data is
$$
U(x,\lambda_1,\lambda_2) = 
\frac{144\lambda_1 e^{-x} + 144\lambda_2 e^{-3x} + 36\lambda_1^2 
\lambda_2 e^{-5x} + 
4\lambda_1 \lambda_2^2 e^{-7x}}{144 + 144\lambda_1^2 e^{-2x} + 
72\lambda_1 \lambda_2 e^{-4x} 
+ 16 \lambda_2^2 e^{-6x} + \lambda_1^2 \lambda_2^2 e^{-8x}}
$$
which for $\lambda_1 = 2, \lambda_2 = 6$ takes the form
$$
U_2(x) = \frac{1}{\cosh x}.
$$

{\bf 4.2.3. The potentials of the Dirac spheres.}

U. Pinkall and J. Richter had constructed the Dirac spheres without using the
the representation theory (\cite{R}). 
We mentioned above
that for these spheres the estimate (\ref{27}) is precise.
Their potentials are
$$
U_N(x) = \frac{N}{2 \cosh x}.
$$
We show above how $U_2(x)$ is obtained via the inverse scattering method 
and it is clear from the trace formula (\ref{21}) and Theorem 7 that 

\begin{proposition}
The discrete spectrum of $U_N(x)$ consists of $\pm \frac{(2j-1)i}{2}$
with $j \leq N$.
\end{proposition}
  
For giving a complete their description of these potentials it needs 
to find the coefficients $\lambda_1,\dots,\lambda_N$. But let us recall that
there exists an infinite family of soliton equations, {\sl the modified
Korteweg--de Vries hierarchy} of nonlinear equations, such that

1) the $m$-th mKdV equation has the form
$$
\frac{\partial}{\partial t_m}U = \frac{\partial^{2m-1}}{\partial x^{2m-1}} U 
+ \dots,
$$
preserves the spectrum of $L$, and transforms $R(k)$ and $\lambda_j$
as follows
$$
\lambda_j \rightarrow \lambda_j \cdot \exp(i2^{2m-1}\kappa_j t_m), 
\ \ \
R(k) \rightarrow R(k) \cdot \exp(i2^{2m-1}k t_m);
$$

2) all flows generated by these equation pairwise commute.

\noindent
Stationary solutions to linear combinations of these flows satisfy
to {\sl the Novikov equations} (\cite{N}):
$$
\left(
a_1 \frac{\partial}{\partial t_1} + \dots + a_m
\frac{\partial}{\partial t_m}
\right) U(x,t_1,\dots) = 0.
$$
These facts are exposed in \cite{AS,FT,Lamb,TS}.
It follows from these formulas that
reflectionless potentials are exactly fast decaying solutions to the 
Novikov equations (for the mKdV hierarchy)
and, given $\kappa_1,\dots,\kappa_N$, the mKdV-orbit of 
a reflectionless potential consists of solutions to these equations
with given $a_1,\dots,a_m$.

\begin{proposition} 
Every reflectionless potential with $\kappa_j = \frac{(2j-1)i}{2}$, where
$ 1 \leq j \leq N$, is obtained from $U_N(x)$ by the mKdV-deformations.
\end{proposition}

{\bf 4.2.4. The Dirac spheres as rational spheres.}

A nice property of the Dirac spheres is that they
are described in terms of rational functions. Indeed, return back to
a plane representation of spheres. It means that we represent $S^2$ as
a complex plane $\C$ completed by a point at infinity.
A conformal parameter $Z$ on $\C$ is related with $x$ and $y$, 
coming in (\ref{23}) and (\ref{31}), as $Z = e^{x+i y}$
and, by (\ref{13}), a spinor field $\Psi(Z,\bar{Z})$ defining via (\ref{4})
a plane representation of a sphere is
$$
\Psi_1(Z,\bar{Z}) = \frac{1}{e^{(x+i y)/2}} \psi_1(x,y), \ \ \
\Psi_2(Z,\bar{Z}) = \frac{1}{e^{(x-i y)/2}} \psi_2(x,y).
$$
Now it is easy to see that $\Psi$ is a rational function of $Z$ and 
$\bar{Z}$.

{\bf 4.2.5. Soliton deformations of reflectionless spheres.}

Each mKdV-flow $U(x,t)$ deforms the eigenfunctions $\varphi^+_1(x,k)$ via
quite simple differential equations which are linear in $\varphi^+_1$.
This gives a deformation of a sphere defined by any linear combination
(\ref{23}). These deformations are called the mKdV-deformations of spheres
and they are reductions (for one-dimensional potentials) of more general
the modified Novikov--Veselov deformations introduced in \cite{Kon}
as deformations of surfaces locally presented via (\ref{4}). 

Here we only mention some interesting facts:

{\sl 1) If $N \geq 2$, then generically a function of the form (\ref{23})
defines not a sphere of revolution and moreover a sphere admitting no
$S^1$-isometries. Nevertheless, soliton deformations of this sphere 
are described by $1+1$-dimensional soliton equations;

2) The Kruskal integrals are first integrals for all mKdV-flows and
for reflectionless spheres they are presented in terms of
$\kappa_1,\dots,\kappa_N$ by the trace formulas (\ref{20});

3) The mKdV-deformations preserve closedness of
spheres (for tori this has been proved in \cite{T2}).}

Generically the values of these first integrals as 
well as the modified Novikov--Veselov deformations depend on a choice of a 
conformal parameter on a sphere and on a surface of revolution 
there exists a distinguished parameter (see \cite{T2}).

{\bf 4.2.6. Reflectionless spheres with} 
$\dim_{\bf H}\mbox{Ker}\, {\cal D} =1$ {\bf and with
large values of the Willmore functional.}

A construction of such spheres is simple: take $\kappa_1 = i/2,
\kappa_2 = -a + it$, and $\kappa_3 = -a + it$, and  fix admissible
$\lambda_j$. Then for each $t > 0$ take a sphere of
revolution $\Sigma_t$ given by $\varphi^+_1(x,i/2)$. By (\ref{21}),
we have
$$
{\cal W}(\Sigma_t) = 4\pi + 8\pi t.
$$
Therefore we conclude

\begin{proposition}
For any $C > 0$ there exists a regular sphere $\Sigma$ in $\R^3$
such that $\dim_{\bf H}\mbox{Ker}\, {\cal D} =1$, $\Sigma$ is defined via
(\ref{4}) by a zero-eigenfunction of ${\cal D}$,  and ${\cal W}(\Sigma) > C$.
\end{proposition} 

We do not check the regularity condition but it is easy using (\ref{31}).

{\bf 4.2.7. Deformations of reflectionless spheres via deformations
of the spectrum.}

We consider a simple example which relates to the previous one (see 4.2.6)
and do that only for demonstrating such deformations.

Take a reflectionless potential $U(x,t$
whose spectral data are 
$\kappa_1,\dots,\kappa_{N+2}$ and
$\lambda_1,\dots,\lambda_{N+2}$ with
$\kappa_{N+1} = -a + i t, \kappa_{N+2} = a + i t$ and
$\lambda_{N+1} = \alpha t, \lambda_{N+2} = \bar{\alpha}t$ with
$a, t \in \R$. 
Take a linear combination (\ref{23}), of eigenfunctions of $L$ 
with eigenvalues of the form $(2m+1)i/2$.
Since $\varphi_1^+(x,\kappa_{N+1})$ and
$\varphi_1^+(x,\kappa_{N+2})$ do not come into this combination, then 
for each $t$ this combination defines a sphere $\Sigma_t$ immersed 
into $\R^3$.

Consider  the limit $t \rightarrow 0$. Then it is easy to see that

{\sl $U(x,t) \rightarrow U(x,0)$ with the spectral data 
$\kappa_1,\dots,\kappa_N$ and $\lambda_1,\dots$,$\lambda_N$, and 
$\Sigma_t$ tends to a sphere $\Sigma_0$ immersed into $\R^3$.}

It is quite sure from the construction that the spectral data of $L$ 
depends on $U(x)$ continuously but we do not know are any estimates for 
stability of the inverse scattering problem for Dirac operators obtained
or not. Hence we may propose the following treating the Dirac spheres as 
``constrained'' Willmore only on a physical rigor level.  

Let $\Sigma$ be a Dirac sphere corresponding to an $N$-soliton potential
and assume that the eigenfunctions corresponding to all levels
of discrete spectrum come into a linear combination (\ref{23}) defining
$\Sigma$. Consider small perturbations of the sphere 
preserving the class of spheres with one-dimensional potentials.
Such perturbations reflect in small perturbations
of the potential and therewith in small perturbations of the discrete
spectrum. These perturbations may result in 
appearance of new eigenvalues and
perturbations of $R(k)$ which transform $U(x)$ into non-reflectionless
potential. But the discrete spectrum $\kappa_j = (2j-1)i/2$ has to be 
preserved because it is all coming in the representation of the sphere.
Now it follows from (\ref{21}) and (\ref{28}) that 
all such perturbations have to increase the value of the Willmore functional. 
By the completeness argument, we conclude that also for all Dirac spheres
and moreover for all spheres such that each eigenvalue of 
$L$ takes the form $(2m+1)i/2$. Therefore we have

{\sl soliton spheres, such that each eigenvalue of $L$ is of the form 
$(2m+1)i/2$,  are critical points of the Willmore functional restricted 
onto the class of spheres with one-dimensional potentials.}

\begin{center}
{\bf \S 5. Final remarks}
\end{center}

1) In the present paper we mostly use a cylindric representation which
enables us to apply well-developed inverse spectral theory for
one-dimensional Dirac operators.

In a plane representation we may treat the
two-dimensional Zakharov--Shabat problem by the $\bar{\partial}$ method
and the nonlocal Riemann problem. In this case we 
consider the spectral data related to one level of energy $E$:
$$
{\cal D} \psi = E\psi,
$$
and in \cite{Bg} the reconstruction problem has been solved for positive $E$
and assuming that the $L_2$-norm of $U$ is sufficiently small.
The last condition is required for a unique solvability 
of the integral equations
which are analogs of (\ref{B2}).
For geometric reasons we have to consider this problem for $E=0$ and 
for potentials with sufficiently large $L_2$-norms, 
$||U||^2_{L_2} \geq \pi$ (the Gauss--Bonnet theorem, (\ref{28})),
and the vanishing of $E$ itself leads to 
appearance of logarithmic singularities of the Green--Faddeev functions 
coming into the kernels of these integral equations.  

In a cylindric representation we may introduce analogs of (\ref{B2})
as follows. 

Find $G(z,k)$ satisfying 
$$
\left[
\left(
\begin{array}{cc}
0 & \partial \\
- \bar{\partial} & 0
\end{array}
\right)
-
\frac{1}{2}
\left(
\begin{array}{cc}
0 & ik \\
ik & 0
\end{array}
\right)
\right]
G(z, k) = \delta(z).
$$
Look for them in the form
$$ 
G(z,k) = 
\left(
\begin{array}{cc}
0 & -e^{-ikx}h_1(z) \\
e^{ikx}h_2(z) & 0
\end{array}
\right)
$$
where 
$\bar{\partial}h_1(z) = \partial h_2(z) = \delta(z)$.
For that consider the Fourier decompositions of $h_1$ and $h_2$:
$$
h_j (z) = 
\int_{-\infty}^{+\infty} d\lambda \sum_{n \in \Z} h^j_{\lambda,n} 
e^{i\lambda x + i n y},
$$
and keeping in mind that
$$
\frac{1}{2\pi}
\int_{-\infty}^{+\infty}
\frac{e^{i\lambda x}}{i\lambda +n} d\lambda =
\cases{
e^{-n x} & for $n>0, x>0$ \cr
-e^{-n x} & for $n<0, x<0$ \cr
0 & for $n<0, x>0$ or $n>0, x<0$},
$$
obtain
$$
h_1(z) = 
\frac{1}{\pi}\left(
\sigma_1(x) +   
\theta(x)\frac{e^{-z}}{1-e^{-z}} + (\theta(x)-1)\frac{e^z}{1-e^z}
\right),
$$
$$
h_2(z) = \frac{1}{\pi}\left(
\sigma_2(x) +   
\theta(x)\frac{e^{-\bar{z}}}{1-e^{-\bar{z}}} + 
(\theta(x)-1)\frac{e^{\bar{z}}}{1-e^{\bar{z}}}
\right).
$$
$$
\partial \sigma_j(x) = \delta(x).
$$
Now choosing $\sigma_j$ as in Appendix B we construct analogs of (\ref{B3})
and (\ref{B4}):
$$
\Phi^-(z,k) = 
\left(
\begin{array}{cc}
0 & e^{-ikx} \\
e^{ikx} & 0
\end{array}\right) -
$$
$$
\frac{1}{\pi}\int_{-\infty}^x d x' \int_0^{2\pi} d y'
\left(
\begin{array}{cc}
0 & \frac{e^{-ik(x-x')}}{e^{z'-z}-1} \\
\frac{e^{ik(x-x')}}{1-e^{\bar{z}'-\bar{z}}} & 0
\end{array}\right)
\cdot U(z') \cdot \Phi^-(z',k) +
$$
$$
\frac{1}{\pi}\int_x^{+\infty}d x' \int_0^{2\pi} d y'
\left(
\begin{array}{cc}
0 & \frac{e^{-ik(x-x')+(z-z')}}{e^{z-z'}-1} \\
\frac{e^{ik(x-x')+(\bar{z}-\bar{z}')}}{1-e^{\bar{z}-\bar{z}'}} & 0
\end{array}\right)
\cdot U(z') \cdot \Phi^-(z',k)
$$
and
$$
\Phi^+(z,k) = 
\left(
\begin{array}{cc}
0 & e^{-ikx} \\
e^{ikx} & 0
\end{array}\right) +
$$
$$
\frac{1}{\pi}\int_x^{+\infty} d x' \int_0^{2\pi} d y'
\left(
\begin{array}{cc}
0 & \frac{e^{-ik(x-x')}}{e^{z-z'}-1} \\
\frac{e^{ik(x-x')}}{1-e^{\bar{z}-\bar{z}'}} & 0
\end{array}\right)
\cdot U(z') \cdot \Phi^+(z',k) -
\label{g10}
$$
$$
\frac{1}{\pi}\int_{-\infty}^x d x' \int_0^{2\pi} d y'
\left(
\begin{array}{cc}
0 & \frac{e^{-ik(x-x')+(z'-z)}}{e^{z'-z}-1} \\
\frac{e^{ik(x-x')+(\bar{z}'-\bar{z})}}{1-e^{\bar{z}'-\bar{z}}} & 0
\end{array}\right)
\cdot U(z') \cdot \Phi^+(z',k).
$$
The kernels of these integral equations
have removable singularities at $z=0$. Indeed,
for $h_1$ the singularity is of the form
$$
\frac{\theta(x)}{z} - \frac{\theta(x)-1}{z} = \frac{1}{z},
$$
i.e., a multiple of a fundamental solution to the Cauchy--Riemann
equations.
For $U$ and $\Phi$ depending on $x$ only, these equations after 
integrating by $y$ reduce to (\ref{B3}) and (\ref{B4}).  

Analogs of (\ref{B5}) for these equations are not quite good because
for sufficiently large $|\mbox{Im}\, k|$ their kernels are exponentially
growing. Nevertheless $\mbox{Ker}\, {\cal D}$ has to admit a description in terms of
solutions to these equations because it is just another representation 
of the spectral problem (\ref{2}) for an operator acting on a
spinor bundle over a two-dimensional sphere.

2) The spinor field $\psi$ on an immersed surface may be obtained also
as the restriction of a parallel spinor field on $\R^3$ onto the
surface. This recently has been exposed in \cite{F} and also 
had been pointed out to the author by Pinkall just after author's talk in 
Amherst (November, 1995) about the paper \cite{T1}
\footnote{Notice that in \cite{T1}
spinors appear as a result of globalizing Weierstrass representations.}. 

Being geometrically invariant this treating of Lemma B does not
enable us to extract spectral-theoretical properties of the Weierstrass
representation which are of a global origin and to develop construction of
surfaces from spectral-theoretical data as we do that in \S 3 and \S 4.   

But this fits into the more general consideration of Dirac operators
on hypersurfaces in $\R^{n+1}$ with induced spin structures. Recently it has
been shown by C. B\"ar (\cite{Baer}) that for 
such hypersurface $M$ there exists at least
$2^{[n/2]}$ distinct eigenvalues $\lambda$ of the induced operator 
${\cal D}$ such that 
$$
\lambda^2 \leq \frac{n^2}{4 \mbox{Vol}\,(M)}{\cal W}(M).
$$ 

3) This work has been done when the author was a guest of SFB 288
in Technische-Universit\"at in Berlin and was also supported by 
INTAS-RFBR (the Russian Foundation for Basic Researches) 
(grant 95-0418).

The author thanks U. Pinkall for helpful discussions.


\vskip10mm

\begin{center}
{\bf Appendix A. Period problem for immersions of surfaces of genus $\geq 2$}
\end{center}

Recall the definition of a global Weierstrass representation of
a surface $\Sigma$ of genus $q \geq 2$ which is conformally equivalent to
$\Sigma_0 = {\cal H}/\Lambda$ with ${\cal H}$ the Lobachevskii 
upper half-plane and $\Lambda$ a lattice in $PSL(2,\R)$ (\cite{T1,T4}):

\begin{definition}
A sphere $\Sigma$ with $g \,(>1)$ handles, immersed into $\R^3$,
possesses a global Weierstrass representation if there exist a real
potential $U$ and functions $\psi_1$ and $\psi_2$, defined on the universal
covering of $\Sigma$, i.e., on ${\cal H}$, such that

1)
\begin{equation}
\left\{
\begin{array}{l}
U(\gamma(z)) = |cz +d|^2 U(z), \\
\psi_1(\gamma(z)) = (cz+d) \psi_1(z), \\
\psi_2(\gamma(z)) = (c\bar{z}+d) \psi_2(z)
\end{array}
\right.
\label{A1}
\end{equation}
for $z \in {\cal H}$ and $\gamma \in \Lambda$, represented by the matrix
$$
\left(
\begin{array}{cc}
a & b \\
c & d
\end{array}
\right),
\ \ \ \
a,b,c,d \in \R, \ \
ad-bc=1;
$$

2) the vector function $\psi$ satisfies (\ref{2}) and 
for a suitable choice of coordinates in $\R^3$ 
defines by (\ref{4}) an immersion of $\Sigma$.
\end{definition}

We have 

\begin{theorem}
Let $\Sigma_0$ be a compact oriented surfaces of genus $g \geq 2$ and let
$U$ and $\psi$ satisfy (\ref{A1}) and define via (\ref{4}) an immersion
of the universal covering of $\Sigma_0$ into $\R^3$. Then this immersion
converts into an immersion of $\Sigma_0$ if and only if 
\begin{equation}
\int_{\Sigma_0} \bar{\psi}_1^2 \, d\bar{z} \wedge \omega =
\int_{\Sigma_0} \psi_2^2 \, d\bar{z} \wedge \omega =
\int_{\Sigma_0} \bar{\psi}_1 \psi_2 \, d\bar{z} \wedge \omega = 0
\label{A2}
\end{equation}
for any holomorphic differential $\omega$ on $\Sigma_0$.
\end{theorem}

A proof of the theorem.

Take a canonical basis 
$\alpha_1,\dots,\alpha_g,\beta_1,\dots,\beta_g$ for $H_1(\Sigma_0)$. This 
means that its intersection form is 
$$
\alpha_j \circ \beta_k = \delta_{j k}, \ \ 
\alpha_j \circ \alpha_k = \beta_j \circ \beta_k = 0
$$
and there exists  loops representing these cycles such that after
cutting along them we obtain a domain $M$ bounded by a polygon 
$$
\partial M =
\alpha_1 \beta_1 \alpha^{-1}_1 \beta^{-1}_1 \dots  
\alpha_g \beta_g \alpha^{-1}_g \beta^{-1}_g.
$$
To this basis corresponds a unique basis of holomorphic differentials
$\omega_1, \dots, \omega_g$ normalized by the condition
$$
\int_{\alpha_k} \omega_j = \delta_{j k}.
$$
Define the period matrix $\Omega$ by
$$
\Omega_{j k} = \int_{\beta_k} \omega_j.
$$
This matrix is symmetric and its imaginary part is positive 
definite (\cite{Fay}). 

It follows from (\ref{4}) and (\ref{A1}) that 
$\bar{\psi}_1^2 d\bar{z}, \psi_2^2 d\bar{z}$, and 
$\bar{\psi}_1 \psi_2 d\bar{z}$ are correctly defined $1$-forms on
$\Sigma_0$ and 
\begin{equation}
\bar{\psi}_1^2 = \bar{\partial} (i X^1 - X^2), \ \ 
\psi_2^2 = \bar{\partial} (i X^1 + X^2), \ \
\bar{\psi}_1 \psi_2 = \bar{\partial} X^3.
\label{A3}
\end{equation}

Introduce the following vectors of translation periods
$$
V_j = \int_{\alpha_j}(\bar{\psi}_1^2 \, d\bar{z} - \bar{\psi}_2^2 d z), \ \
V_{g+j} = \int_{\beta_j}(\bar{\psi}_1^2 \, d\bar{z} - \bar{\psi}_2^2 d z)
$$
and
$$
W_j =  \int_{\alpha_j}(\psi_1 \bar{\psi}_2 \, d z + 
\bar{\psi}_1 \psi_2 \, d\bar{z}), \ \ 
W_{g+j} =  \int_{\beta_j}(\psi_1 \bar{\psi}_2 \, d z + 
\bar{\psi}_1 \psi_2 \, d\bar{z}).  
$$ 

Denote by $Y$ and $Z$ the vectors
$$
Y_j =\int_{\Sigma_0} \bar{\psi}_1^2 \, d\bar{z} \wedge \omega_j , \ \
Y_{j+g} = -
\overline{\int_{\Sigma_0} \psi_2^2 \, d\bar{z} \wedge \omega_j}
$$
and
$$
Z_j = 
\int_{\Sigma_0} \bar{\psi}_1 \psi_2 \, d\bar{z} \wedge \omega_j,\ \
Z_{j+g} = 
\overline{\int_{\Sigma_0} \bar{\psi}_1 \psi_2 \, d\bar{z} \wedge \omega_j}.
$$

It is evident that an immersion of ${\cal H}$ is converted into an 
immersion of $\Sigma_0$ if and only if $V = W = 0$.

Now, by the Stokes theorem and (\ref{A3}), we have
$$
\int_{\Sigma_0} \bar{\psi}_1^2 \, d\bar{z} \wedge \omega_j = 
\int_{\partial M} (i X^1 - X^2)\omega_j = 
\sum_{k=1}^g
\left( V_k \int_{\beta_k} \omega_j - V_{k+g}\int_{\alpha_k} \omega_j\right),
$$
\begin{equation}
\int_{\Sigma_0} \psi_2^2 \, d\bar{z} \wedge \omega_j = 
\int_{\partial M} (i X^1 + X^2)\omega_j = 
\sum_{k=1}^g
\left( -\bar{V}_k \int_{\beta_k} \omega_j + 
\bar{V}_{k+g}\int_{\alpha_k} \omega_j\right)
\label{A4}
\end{equation}
$$
\int_{\Sigma_0} \bar{\psi}_1 \psi_2 \, d\bar{z} \wedge \omega_j = 
\int_{\partial M} X^3\omega_j = 
\sum_{k=1}^g
\left( W_k \int_{\beta_k} \omega_j - 
W_{k+g}\int_{\alpha_k} \omega_j\right).
$$
Consider the $2g \times 2g$-matrix
$$
\tilde{\Omega} = 
\left(
\begin{array}{cc}
\Omega & -1\\
\bar{\Omega} & -1
\end{array}
\right).
$$
Since $\mbox{Im}\, \Omega$ is positive definite, $\tilde{\Omega}$ is 
nondegenerate. Rewrite (\ref{A4}) as follows 
$$
\tilde{\Omega} V = Y,\ \ \ \tilde{\Omega} W = Z,
$$
and conclude that $V = W =0$ if and only if 
\begin{equation}
Y = Z = 0.
\label{A5}
\end{equation}
Since $\omega_j$ form a basis for holomorphic differentials,
(\ref{A5}) is equivalent to vanishing integrals (\ref{A2}) for any
holomorphic differential on $\Sigma_0$.

This proves the theorem.

Proposition 4 from \cite{T4} which settles the period
problem for tori may be reformulated as (\ref{A2}).

\vskip10mm

\begin{center}
{\bf Appendix B. The inverse scattering problem for the one-dimensional
Dirac operator}
\end{center}

{\bf B1. The forward scattering problem.}

Consider the linear problem
\begin{equation}
L \varphi = 
\left(
\begin{array}{cc}
0 & ik \\
ik & 0
\end{array}
\right)
\varphi
\label{B1}
\end{equation}
with
$$
L = 
\left(
\begin{array}{cc}
0 & \partial_x \\
- \partial_x & 0
\end{array}
\right)
+
2
\left(
\begin{array}{cc}
U & 0 \\
0 & U
\end{array}
\right)
$$
and $k$ a spectral parameter.
This is the simplest reduction of the Zakharov--Shabat linear problem
corresponding to the case when both potentials are equal to $U(x)$ and are
real-valued.

For each $k \in \R \setminus\{0\}$ 
the system (\ref{B1}) has a two-dimensional space of
solutions. Take a matrix $\Phi(x,k)$ whose columns form a basis for this space.
For $U = 0$ we take
$$
\Phi_0(x,k) = 
\left(
\begin{array}{cc}
0 & e^{-ikx} \\
e^{ikx} & 0
\end{array}\right)
$$
For fast decaying $U$ 
the matrix
$\Phi$ may be defined by the integral equation
\begin{equation}
\Phi(x,k) = \Phi_0(x,k) - \int_{-\infty}^{+\infty} G(x-x',k) \cdot 2U(x') 
\cdot \Phi(x',k) \, d x'
\label{B2}
\end{equation}
where
$$
\left[
\left(
\begin{array}{cc}
0 & \partial_x \\
- \partial_x & 0
\end{array}
\right)
-
\left(
\begin{array}{cc}
0 & ik \\
ik & 0
\end{array}
\right)
\right]G(x,k) = \delta(x),
$$
i.e., $G(x,k)$ is a fundamental solution to 
(\ref{B1}) with $U = 0$  
and it may be taken in the form
$$
G(x,k) = 
\left(
\begin{array}{cc}
0 & -e^{-i k x} g_1(x) \\
e^{i k x} g_2(x) & 0
\end{array}
\right)
$$
with
$$
\partial_x g_j(x) = \delta(x).
$$
The last equality is satisfied exactly by 
$\theta(x) + \mbox{const}$
with $\theta(x)$ the Heaviside function
$$
\theta(x) = 
\cases{ 
1& for $x \geq 0$ \cr
0& for $x < 0$}.
$$
This freedom to choose $g_j$ enables us
to define solutions to (\ref{B1}) converging to the free waves, given
by the columns of $\Phi_0(x,k)$,
as $x \rightarrow -\infty$ or $x \rightarrow +\infty$.
Indeed, for  $g_1(x) = g_2(x) = \theta(x)$ the system (\ref{B2})
takes the form
\begin{equation}
\Phi^- (x,k) =
\left(
\begin{array}{cc}
0 & e^{-i k x} \\
e^{i k x} & 0
\end{array}\right)
-
\label{B3}
\end{equation}
$$
\int_{-\infty}^x
\left(
\begin{array}{cc}
0 & -e^{-i k (x-x')} \\
e^{i k (x-x')} & 0
\end{array}\right)
\cdot 2U(x') \cdot \Phi^-(x',k)
\, d x'
$$
and for $g_1(x) = g_2(x) = \theta(x) -1$ it is
\begin{equation}
\Phi^+ (x,k) =
\left(
\begin{array}{cc}
0 & e^{-i k x} \\
e^{i k x} & 0
\end{array}\right)
+
\label{B4}
\end{equation}
$$
\int_x^{+\infty}
\left(
\begin{array}{cc}
0 & -e^{-i k (x-x')} \\
e^{i k (x-x')} & 0
\end{array}\right)
\cdot 2U(x') \cdot \Phi^+(x',k)
\, d x'.
$$
The solutions to (\ref{B3}) and (\ref{B4}) are called the Jost
functions. 

1) {\sl For each $k \in \R \setminus \{0\}$ the equations (\ref{B3}) and 
(\ref{B4}) have unique solutions.}

These equations have the form 
$$
\Phi^{\pm}(x,k) = \Phi_0(x,k) + (A^{\pm}\circ \Phi^{\pm})(x,k)
$$
where operators $A^{\pm}$ are of the Volterra type. Solutions to them are
given by the Neumann series
$$
\Phi^{\pm}(x,k) = \sum_{j=0}^{\infty} \Phi^{\pm}_j (x,k),
$$
with
$\Phi^{\pm}_j(x,k) = (A^{\pm} \circ \Phi^{\pm}_{j-1})(x,k)$,
which uniformly converge in $x$ on each compact interval.

Denote the $j$-th column of $\Phi^{\pm}(x,k)$ by $\varphi^{\pm}_j$.
Each of the pairs $(\varphi^-_1,\varphi^-_2)$ and 
$(\varphi^+_1,\varphi^+_2)$ forms a basis for 
solutions to (\ref{B1}).

2) {\sl The functions $\varphi^-_1(x,k)e^{-ikx}$ and $\varphi^+_2(x,k)e^{ikx}$
are analytically continued onto the lower-half plane $\mbox{Im}\, k <0$
and the functions  $\varphi^-_2(x,k)e^{ikx}$ and $\varphi^+_1(x,k)e^{-ikx}$
are analytically continued onto the upper-half plane $\mbox{Im}\, k >0$.}

Without loss of generality we explain this fact for $f(x,k) = 
\varphi^-_1(x,k)e^{-ikx}$. This function satisfies the integral equation
\begin{equation}
f(x,k) = 
\left(
\begin{array}{c}
0 \\ 1
\end{array}\right)
-
\int^x_{-\infty}
\left(
\begin{array}{cc}
0 & - e^{-2ik(x-x')} \\
1 & 0
\end{array}\right)
2U(x')f(x',k) \, d x'
\label{B5}
\end{equation}
For $\mbox{Im}\, k < 0$ its kernel decays exponentially as $x \rightarrow 
\infty$ and the Neumann series for this equation converge.
 
3) {\sl Given vector functions $\theta(x) = (\theta_1(x), \theta_2(x))$
and $\tau(x) = (\tau_1(x), \tau_2(x))$ satisfying (\ref{B1}),
the Wronskian $W(\theta,\tau)(x) = \theta_1(x)\tau_2(x) - 
\theta_2(x)\tau_1(x)$ is constant, i.e., independent of $x$.}

This Wronskian identity is obtained by straightforward computations
and implies that 
$$
\det \Phi^{\pm}(x,k) = -1.
$$

4) {\sl For $k \in \R \setminus \{0\}$ the matrices 
$\Phi^+(x,k)$ and $\Phi^-(x,k)$ are related as
$$
\Phi^-(x,k) = S(k)\Phi^+(x,k)
$$
with the scattering matrix $S(k)$ independent of $x$ with} 
$$
\det S(k) = 1.
$$

Indeed, the columns of $\Phi^+(x,k)$ and $\Phi^-(x,k)$ form different bases
for solutions to (\ref{B1})
and, therefore, are linearly dependent:
$$
\varphi^-_1(x,k) = s_{11}(k)\varphi^+_1(x,k) + s_{12}(k)\varphi^+_2(x,k),
$$
$$  
\varphi^-_2(x,k) = s_{21}(k)\varphi^+_1(x,k) + s_{22}(k)\varphi^+_2(x,k).
$$
Since $\det \Phi^{\pm}(x,k) = -1$, we have $\det S(k) = 1$.

Denote $s_{22}(k)$ by $a(k)$ and $s_{21}(k)$ by $b(k)$.

Notice that, if  
$\theta(x) = (\theta_1(x),\theta_2(x))$ satisfies (\ref{B1}), then 
the function $\tilde{\theta}(x) = 
(\theta_2(x), - \theta_1(x))$ satisfies (\ref{B1}) with $-k$ substituted for
$k$. This implies that

5) {\sl For $k \in \R \setminus \{0\}$,}
\begin{equation}
\Phi^{\pm}(x,k) = 
J \cdot \Phi^{\pm}(x,-k) \cdot J 
\ \ \ 
\mbox{with} \ \ 
J = 
\left(
\begin{array}{cc}
0 & 1 \\ -1 & 0
\end{array}
\right).
\label{B6}
\end{equation} 

Since $U(x)$ is real-valued, it is also clear that
\begin{equation}
\Phi^{\pm}(x,k) = \overline{\Phi^{\pm}(x,-k)} \ \ \ 
\mbox{for $k \in \R \setminus \{0\}$}.
\label{B7}
\end{equation}

6) {\sl The scattering matrix takes the form}
$$
S(k) = 
\left(
\begin{array}{cc}
\overline{a(k)} & -\overline{b(k)} \\
b(k) & a(k)
\end{array}
\right)
\ \ \ \
\mbox{with
$|a(k)|^2 + |b(k)|^2 = 1$}.
$$

Notice that (\ref{B6}) implies
$\Phi^-(x,-k) = \left( -J S(k)J\right) \Phi^+(x,-k)$
and (\ref{B7}) implies
$\Phi^-(x,-k) = \overline{S(k)} \Phi^+(x,-k)$ 
It follows from these equalities that
$\overline{S(k)} = -J S(k)J$ which proves 6).
The following quantities
$$
T(k) = \frac{1}{a(k)}, \ \ \ 
R(k) = \frac{b(k)}{a(k)}
$$
are called {\sl the transmission coefficient} and 
{\sl the reflection coefficient} respectively.
It is shown that $a(k)$ vanishes nowhere on $\R \setminus \{0\}$. 

7) {\sl $T(k)$ is analytically continued onto
the upper-half plane $\mbox{Im}\, k \geq 0$.}

Indeed, this follows from 2) and 
\begin{equation}
a(k) = W(\varphi^+_1(x,k),\varphi^-_2(x,k)) =
W(e^{-i k x}\varphi^+_1(x,k),e^{i k x}\varphi^-_2(x,k)).
\label{B8}
\end{equation}

The poles of $T(k)$ correspond to {\sl bounded states},
i.e., to solutions to (\ref{B1}) which decay exponentially as
$x \rightarrow \pm \infty$. These solutions are  
$\varphi^+_1(x,\kappa)$ or $\varphi^-_2(x,\kappa)$ where
$\kappa$ is a pole of $T(k)$ and, since $a(\kappa) = 0$,
these functions are linearly dependent
\begin{equation}
\varphi^-_2(x,\kappa) = \mu(\kappa) \varphi^+_1(x,\kappa),
\ \ \ \mu(\kappa) \in \C.
\label{B9}
\end{equation}  
Some computations lead to the conclusion which
we only recall:

8) {\sl $T(k)$ has only simple poles in $\mbox{Im}\, k > 0$ and for fast 
decaying, for instance, for exponentially decaying, potentials there are
finitely many poles of $T(k)$.}

The transform
$*: (\xi_1,\xi_2) \rightarrow
(\bar{\xi}_2, -\bar{\xi}_1)$
maps $\varphi^{\pm}_1(x,k)$ into a multiple of $\varphi^{\pm}_2(x,\bar{k})$.
Therefore either $\varphi^+_1(x,\kappa)$ and $\varphi^+_2(x,\bar{\kappa})$ 
both decay exponentially as $x \rightarrow \pm \infty$ or neither of them do.
This implies that 

9) {\sl The discrete spectrum of $L$ is 
preserved by the complex conjugation $\kappa \rightarrow \bar{\kappa}$.}

The following quantities form {\sl the spectral data} of $L$.

{\bf Spectral data:}

{\sl 
1) the reflection coefficient $R(k)$ ;

2) the poles of $T(k)$: $\kappa_1,\dots,\kappa_N$ ;

3) the products  $\lambda_j=i\gamma_j \mu_j$, where
$\gamma_j = \gamma(\kappa_j)$ are 
the residues of $T(k)$ at $\kappa_j$
and $\mu_j = \mu(\kappa_j)$ relate $\varphi^+_1$
and $\varphi^-_2$ at $k = \kappa_j$ (\ref{B9}).}

Notice, that $L$ has a continuous spectrum $k \in \R \setminus 
\{0\}$ of multiplicity two and a discrete spectrum 
$\kappa_1,\dots,\kappa_N,\bar{\kappa}_1,\dots,\bar{\kappa}_N$ of 
multiplicity one.

Since $U(x)$ is real-valued, this reflects
in ``reality conditions'' met by the spectral data.

First, notice that 
\begin{equation}
\varphi^{\pm}_j(x,-k) = \overline{\varphi^{\pm}_j(x,k)}
\ \ \
\mbox{for $k \in \R \setminus \{0\}$}.
\label{B10}
\end{equation}
It follows from (\ref{B8}) and (\ref{B10}) that 
$a(k) = \overline{a(-k)}$ for $k \in \R \setminus \{0\}$}.
Consider now the meromorphic function
$F(k) = T(k) - \overline{T(-\bar{k})}$ 
defined on $\{\mbox{Im}\, k \geq 0\}$.
We see that it vanishes everywhere on the boundary of 
$\{\mbox{Im}\, k > 0\}$ and, therefore, $F(k) \equiv 0$.
This implies that 

{\sl 10) The poles of $T(k)$ are symmetric with respect to the imaginary 
axis , i.e., if $(k_R + i k_I)$ is a pole of $T(k)$, where
$k_R, k_I \in \R$, then $(-k_R + i k_I)$ is also the pole of $T(k)$.
Moreover the residues of $T(k)$ are related as 
$\gamma(k_R+i k_I) = -\overline{\gamma(-k_R+i k_I)}$
and if $\Re \kappa_j = 0$ then 
$\mbox{Re}\, \gamma(\kappa_j)=0$.} 

Consider the analytic continuation
of the functions 
$G_1(x,k) = \varphi^+_1(x,k) - \overline{\varphi^+_1(x,-\bar{k})}$
and 
$G_2(x,k) = \varphi^-_2(x,k) - \overline{\varphi^+_2(x,-\bar{k})}$
onto the same upper-half plane. Since, by (\ref{B10}), they identically
vanish on the boundary, we conclude that
$G_1(x,k) \equiv G_2(x,k) \equiv 0$.
This implies that

{\sl 11) The coefficients $\mu_(\kappa)$ relating, via (\ref{B9}),
the functions $\varphi^+_1(x,k)$ and  $\varphi^-_2(x,k)$ 
at the poles of $T(k)$ are complex conjugate, i.e. 
$\mu(k_R + i k_I) = \overline{\mu(-k_R + i k_I)}$,
and if $\mbox{Re}\, \kappa_j = 0$ then $\mu(\kappa_j) \in \R$.}

We summarize 10) and 11) in 

{\sl R1) The poles of $T(k)$ are symmetric with respect to the imaginary axis,
$\lambda_j = \bar{\lambda}_l$ for a pair of symmetric poles $\kappa_j$ and 
$\kappa_l$, and if $\mbox{Re}\, \kappa_m =0$ then
$\lambda_m \in \R$.}  

The definition of $T(k)$ and (\ref{B10}) also imply that

{\sl R2) $R(k) = \overline{R(-k)}$.}    

The conditions R1--R2 perfectly distinguish the spectral 
data of real potentials. 

{\bf B2. The Marchenko equations.}

Introduce the following representations for
$\varphi^-_2(x,k)$ and $\varphi^+_1(x,k)$
\footnote{The existence of these representations is established via
the Goursat equations which we discuss in B3.}
:
$$
\varphi^-_2(x,k) = 
\left(
\begin{array}{c}
e^{-ikx} \\ 0
\end{array}
\right)
+
\int_{-\infty}^x d x'
\left(
\begin{array}{c}
A_1(x,x') \\ A_2(x,x')
\end{array}
\right)
e^{-i k x'}
$$
and
\begin{equation}
\varphi^+_1(x,k) = 
\left(
\begin{array}{c}
0 \\ e^{ikx}
\end{array}
\right)
+
\int_x^{+\infty} d x'
\left(
\begin{array}{c}
B_1(x,x') \\ B_2(x,x')
\end{array}
\right)
e^{i k x'}.
\label{B11}
\end{equation}
It follows from (\ref{B6}) that
\begin{equation}
\varphi^+_2(x,k) = 
\left(
\begin{array}{c}
e^{-ikx} \\ 0
\end{array}
\right)
+
\int_x^{+\infty} d x'
\left(
\begin{array}{c}
B_2(x,x') \\ -B_1(x,x')
\end{array}
\right)
e^{-i k x'}.
\label{B12}
\end{equation}
By the definition of $T(k)$ and $R(k)$, we have
\begin{equation}
T(k) \varphi^-_2(x,k) = R(k) \varphi^+_1(x,k) + \varphi^+_2(x,k).
\label{B13}
\end{equation}
Introduce the Fourier transforms 
of $R(k)$ and $(T(k)-1)$:
$$
\Gamma(z) = \frac{1}{2\pi}
\int_{-\infty}^{+\infty}
(T(k)-1)e^{-i k z} \, d k,
\ \ \ \ 
r(z) = 
\frac{1}{2\pi}\int_{-\infty}^{+\infty}
R(k)e^{-i k z} \, d k,
$$
and also apply the Fourier transformation to both sides of (\ref{B13}), 
obtaining a pair of integral equations corresponding to 
the coefficients of $\varphi$:
\begin{equation}
\Gamma(x+t) + \int_{-\infty}^x \Gamma(x'+t)A_1(x,x') \, d x' + 
\theta(x+t)A_1(x,-t) = 
\label{B14}
\end{equation}
$$
\theta(-t-x)B_2(x,-t) + \int_x^{+\infty}r(x'-t)B_1(x,x')\, d x'.
$$
and
\begin{equation}
\int_{-\infty}^x \Gamma(x'+t)A_2(x,x') \, d x' + 
\theta(x+t)A_2(x,-t) =
\label{B15}
\end{equation}
$$
r(x-t)-\theta(-t-x)B_1(x,-t) + \int_x^{+\infty}r(x'-t)B_2(x,x')\, d x'.
$$

Since $e^{-i k x}\varphi^+_1(x,k)$ and $e^{i k x}\varphi^-_2(x,k)$ 
converge very fast to 
$\left(\begin{array}{c} 0 \\ 1 \end{array}\right)$ 
and 
$\left(\begin{array}{c} 1 \\ 0  \end{array}\right)$, respectively, 
as $k \rightarrow \infty$, Im $k \geq 0 $, 
we conclude from (\ref{B8}) that $(T(k) -1)$ is fast decaying and for 
$z < 0$ we have
$$
\Gamma(z) = \frac{1}{2\pi}
\int_{-\infty}^{+\infty}
(T(k)-1)e^{-i k z}\, d k = \sum_j i\gamma_j e^{-i\kappa_j z}
$$
where $\gamma_j$ is the residue of $T(k)$ at $\kappa_j$.
If there are no poles of $T(k)$, then $\Gamma(z) = 0$ for $z<0$ and
for $x+t <0$ the left-hand sides of  (\ref{B14}) and (\ref{B15}) 
vanish. Otherwise
for $x+t < 0$ the Fourier transform of
the left-hand side of (\ref{B13}) is
$$
\sum_{j=1}^N i\gamma_j e^{-i\kappa_j t} \varphi^-_2(x,\kappa_j).
$$
Since $a(\kappa_j)=0$, we have (see (\ref{B9}))
$$
\varphi^-_2(x,\kappa_j) = \mu_j \varphi^+_1(x,\kappa_j), \ \ \ \mu_j \in \C,
$$
and substituting (\ref{B11}) into (\ref{B14}) and (\ref{B15}) we obtain
$$
\int_x^{+\infty}
\sum_{j=1}^N i\gamma_j e^{-i\kappa_j t} 
\mu_j B_1(x,x')e^{i\kappa_j x'}\, d x' =
$$
$$
\theta(-t-x)B_2(x,-t) + \int_x^{+\infty}r(x'-t)B_1(x,x')\, d x'
$$
and
$$
\sum_{j=1}^N i\gamma_j e^{-i\kappa_j t} 
\mu_j e^{i\kappa_j x} +
\int_x^{+\infty}
\sum_{j=1}^N i\gamma_j e^{-\kappa_j t} 
\mu_j B_2(x,x')e^{i\kappa_j x'} \, d x' =
$$
$$ 
r(x-t)-\theta(-t-x)B_1(x,-t) + \int_x^{+\infty}r(x'-t)B_2(x,x')\, d x'.
$$
where $x+t <0$.
Introducing the function
$$
\Omega(z) = r(z) - \sum_{j=1}^N i\gamma_j\mu_j e^{i\kappa_j z} =
r(z) -  \sum_{j=1}^N \lambda_j e^{i\kappa_j z}
$$
and substituting $-y$ for $t$, 
rewrite these equations as follows
\begin{equation}
\cases{
B_2(x,y) + \int_x^{+\infty} B_1(x,x')\Omega(x'+y) \, d x' = 0, \cr
\Omega(x+y) - B_1(x,y) + \int_x^{+\infty} B_2(x,x') \Omega(x'+y)\, d x' = 0
}
\label{B16}
\end{equation}
where $y > x$. 
These equations 
are called {\sl the Marchenko equations} (for the Za\-kha\-rov--Sha\-bat
linear problem).

Notice that $\Omega(x)$ is uniquely reconstructed from the 
spectral data of $L$ and if the 
reality conditions R1--R2 hold then $\Omega(x)$ is real-valued
for $x \in \R$.

{\bf B3. The inverse scattering problem.}

Assume that $\Omega(z)$ is constructed from the spectral data
of $L$ and the equations (\ref{B16})
are solved, i.e., the functions $B_1(x,y)$
and $B_2(x,y)$ are known for $y > x$. 
In fact, these equations are of the Volterra type and uniquely solvable.
Moreover 
\begin{equation}
\lim_{y \rightarrow +\infty} B_1(x,y) = 
\lim_{y \rightarrow +\infty} B_2(x,y) = 0,
\label{B17}
\end{equation}
the limits of $B_j(x,y)$ as $y \rightarrow x$ are
defined and we denote them by $B_j(x,x)$.

Substituting (\ref{B12}) into (\ref{B1}),
we obtain a pair of equations corresponding to the rows of a matrix
equation. The first of them is
$$
- \int_x^{+\infty} \frac{B_1(x,y)}{\partial x}e^{-i k y}\, d y + 
e^{-i k x}B_1(x,x) + 
i k\int_x^{+\infty} B_1(x,y)e^{-i k y}\, d y +
$$
$$
2U(x)e^{-i k x} + 2U(x) \int_x^{+\infty} B_2(x,y)e^{-i k y}\, d y = 0.
$$
Integrating by parts
$$
i k\int_x^{+\infty} B_1(x,y)e^{-i k y}\, d y = 
$$
$$
\int_x^{+\infty} \frac{\partial B_1(x,y)}{\partial y}e^{-i k y}\, d y
- \lim_{y \rightarrow +\infty}(B_1(x,y)e^{-i k y}) + B_1(x,x)e^{-i k x},
$$
we finally derive
$$
\int_x^{+\infty} 
\left(\frac{\partial B_1(x,y)}{\partial y} -
\frac{\partial B_1(x,y)}{\partial x} + 2U(x)B_2(x,y) \right)e^{-i k y}\, d y +
$$
$$
2(U(x)+B_1(x,x))e^{-i k x} = 0
$$
which implies
\begin{equation}
\frac{\partial B_1(x,y)}{\partial y} -
\frac{\partial B_1(x,y)}{\partial x} + 2U(x)B_2(x,y) = 0
\label{B18}
\end{equation}
and
\begin{equation} 
U(x) = - B_1(x,x).
\label{B19}
\end{equation}
Analogously we infer that the second equation is equivalent to
\begin{equation}
\frac{\partial B_2(x,y)}{\partial x} + 
\frac{\partial B_2(x,y)}{\partial y}  = -2U(x)B_1(x,y).
\label{B20}
\end{equation}
Substituting $x$ for $y$ in (\ref{B20}) and taking (\ref{B19})
in account, we derive
\begin{equation} 
\frac{d B_2(x,x)}{d x} = 2U^2(x).
\label{B21}
\end{equation}

In fact, 
the integral representation (\ref{B12}) is initially derived from 
the Goursat equations  
(\ref{B18}) and (\ref{B20}) with the boundary conditions 
(\ref{B17}) and (\ref{B19}).

The formula (\ref{B19}) gives a solution to the inverse scattering 
problem: reconstructing the potential from the spectral data.
If the conditions R1--R2 hold 
then $\Omega(x)$ is real-valued for
$x \in R$ and the solution to (\ref{B16}) is also real-valued for
$x, y \in \R$.

To complete this scheme it needs to prove that starting from the spectral data
we obtain via (\ref{B16}) and (\ref{B19}) a potential $U(x)$ with the same
spectral data. The main tool in proving that is the fact that (\ref{B16})
is just the Fourier transform of (\ref{B13}). A detailed analysis would
enable us to distinguish the decay of a potential in terms of the data.
For the one-dimensional Schr\"odinger operator on the line, i.e., the most
similar problem to (\ref{B1}), this had been done in \cite{F1}
(see, also \cite{M}) where one can find a detailed study of this problem.

\newpage

\end{document}